\documentclass[11pt]{article}
\usepackage{amssymb}

\begin{document}

\begin{center}
\textbf{MEAN CURVATURE MOTION OF GRAPHS WITH CONSTANT CONTACT ANGLE AND MOVING BOUNDARIES}\\
by A. Freire\\

University of Tennessee, Knoxville \end{center}

\begin{abstract}
We consider the motion by mean curvature of an $n$-dimensional graph
over a time-dependent domain in $\mathbb{R}^n$, intersecting
$\mathbb{R}^n$ at a constant angle. In the general case, we prove
local existence for the corresponding quasilinear parabolic equation
with a free boundary, and derive a continuation criterion based on
the second fundamental form. If the initial graph is concave, we
show this is preserved, and that the solution exists only for finite
time. This corresponds to a symmetric version of mean curvature
motion of a network of hypersurfaces with triple junctions, with
constant contact angle at the junctions.
\end{abstract}

\vspace{.5cm}

\textbf{1. Time-dependent graphs with a contact angle condition.}
\vspace{.2cm}

We consider a moving hypersurface $\Sigma_t$ in $\mathbb{R}^{n+1}$,
with normal velocity equal to its mean curvature, assumed to be a
graph over a time-dependent open set $D(t)\subset {\mathbb R}^n$
(not necessarily bounded, or connected.) The (properly embedded)
$(n-1)$-submanifold of intersection:
$$\Gamma(t)=\Sigma_t\cap {\mathbb R^n}=\partial D(t)$$
is a `moving boundary'. Along $\Gamma(t)$ we impose a constant-angle
condition:
$$\langle N,e_{n+1}\rangle_{|\Gamma(t)}=\beta,$$
where $0<\beta<1$ is a \emph{constant} and $N$ is the upward unit
normal of $\Sigma_t$. `Mean curvature motion' is defined by the law:
$$V_N=H,$$
where $V_N=\langle V,N\rangle$, with $V=\partial_tF$ the velocity
vector in a given parametrization $F(t)$ of $\Sigma_t$ ($V$ depends
on the parametrization, while $V_N$ does not). A particular
parametrization yields `mean curvature flow':
$$\partial_tF=HN.$$
For graphs, it is natural to consider `graph mean curvature motion':
if $\Sigma_t=\mbox {graph }w(t)$ for a function
$w(t):D(t)\rightarrow {\mathbb R}$, imposing $\langle
\partial_tF,N\rangle=H$ with $F(y,t)=[y,w(y,t)]$ for $y\in D(t)$, we
find:
$$w_t=\sqrt{1+|D w|^2}H$$
(and the velocity is vertical, $\partial_tF=w_te_{n+1})$. With the
contact angle condition, we obtain a free boundary problem for a
quasilinear PDE:
\newpage

$$\left\{ \begin{array}{l}w_t=g^{ij}(Dw)w_{ij}\quad \mbox{ in } D(t),\\
w=0,\quad \beta\sqrt{1+|Dw|^2}=1\mbox{ on }\partial
D(t),\end{array}\right .$$ where $g^{ij}
(Dw)=\delta^{ij}-w_iw_j/(1+|Dw|^2)$ is the inverse metric
matrix.\vspace{.2cm}

\emph{Remark 1.1.} It is easy to see that the constant-angle
boundary condition is \emph{incompatible} with mean curvature flow
parametrized over a \emph{fixed} domain $D_0$: on $\partial D_0$ we
would have $\langle F,e_{n+1}\rangle=0$, leading to $\langle
\partial_tF,e_{n+1}\rangle=0$, incompatible with $\partial_tF=HN$ and
$\langle N,e_{n+1}\rangle=\beta$. If we parametrize over
time-dependent domains, mean curvature flow and graph m.c.m. lead to
identical normal velocities for the moving boundary (see section
2.)\vspace{.2cm}

 To establish short-time existence (in parabolic H\"{o}lder
spaces) we will work with a third realization of the motion, defined
over a fixed domain:
$$F(t):D_0\rightarrow {\mathbb R}^{n+1},\quad
F(x,t)=[\varphi(x,t),u(x,t)]\in {\mathbb R}^n\times {\mathbb R},$$
where $\varphi(t):D_0\rightarrow D(t)$ is a diffeomorphism and $F$
is a solution of the parabolic system:
$$F_t=g^{ij}(DF)F_{ij},$$
where $g_{ij}=\delta_{ij}+\langle F_i,F_j\rangle$ is the induced
metric on $\Sigma_t$ and $g^{ij}$ is the inverse metric matrix.
\vspace{.2cm}

In the first part of the paper (sections 3 to 8) we prove the
following short-time existence theorem (on $Q:=D_0\times
[0,T]$:\vspace{.2cm}

\textbf{Theorem 1.1.} Let $\Sigma_0\subset {\mathbb R}^{n+1}$ be a
$C^{3+\alpha}$ graph over $D_0\subset {\mathbb R}^n$ satisfying the
contact and angle conditions at $\partial D_0$. There exists a
parametrization $F_0=[\varphi_0,u_0]\in C^{2+\alpha}(D_0)$ of
$\Sigma_0$, $T>0$ depending only on $F_0$ and a unique solution
$F\in C^{2+\alpha,1+\alpha/2}(Q^T;{\mathbb R}^{n+1})$ of the system:
$$\left \{ \begin{array}{l}
\partial_t F-g^{ij}(DF)\partial_i\partial_jF=0,\\
u_{|\partial D_0}=0, \quad N^{n+1}(D\varphi,Du)_{|\partial
D_0}=\beta,\end{array}\right .$$ with initial data $F_0$ and
satisfying, in addition, the `orthogonality conditions' at $\partial
D_0$ (described in section 3.) \vspace{.2cm}

The system and boundary conditions are discussed in more detail in
section 3. Sections 4, 5, and 6 deal with compatibility at $t=0$,
linearization and the verification that the boundary conditions
satisfy `complementarity'. In particular, adjusting the initial
diffeomorphism $\varphi_0$ to ensure compatibility (section 4) leads
to the `loss of differentiability' seen in theorem 1.1. The required
estimates in H\"{o}lder spaces for the linearized system are
described in section 7, and the proof concluded (by a fixed-point
argument) in section 8. While the general scheme is standard, due to
the particular boundary conditions adopted many details had to be
worked out from first principles. Free boundary-type problems for
mean curvature motion of graphs have apparently not previously been
considered.

We describe the evolution equations in the rotationally symmetric
case in section 9 (including a stationary example for the exterior
problem) and the extension to the case of a graph motion $\Sigma_t$
intersecting fixed support hypersurfaces orthogonally in section 10.
\vspace{.2cm}

The original motivation for this work was to establish (by classical
parabolic PDE methods) existence-uniqueness for mean curvature
motion of networks of surfaces meeting along triple junctions with
constant-angle conditions. One can use a motion $\Sigma_t$ of graphs
with constant contact angle to produce examples of `triple junction
motion': three hypersurfaces moving by mean curvature meeting along
an ($n-1$)-dimensional submanifold $\Sigma(t)$ so that the three
normals make constant angles (say, 120 degrees) along $\Gamma(t)$.
The simplest way to do this is by \emph{reflection}  on $\mathbb
R^n$, so the hypersurfaces are $\Sigma_t, \bar{\Sigma_t}$ and
$\mathbb{R}^n-\bar{D}(t)$. If $\Sigma_t=\mbox{graph }w(t)$ with
$w>0$, the system is \emph{embedded} in ${\mathbb R}^{n+1}$. This is
mean curvature motion of a `symmetric triple junction of
graphs'.\vspace{.2cm}

Short-time existence holds for general triple junctions of graphs
moving by mean curvature with constant 120 degree angles at the
junction, provided a compatibility condition holds along the
junction (see section 15). The idea of proof is similar to the one
given here; since the details are easier to understand in the
symmetric case, we decided to do this first. In addition, in the
present case it is possible to go a lot further towards a geometric
global existence result. Motivated by recent work on `lens-type'
curve networks \cite{LensSeminar}, in the second part of the paper
(sections 11-14) we consider continuation criteria and preservation
of concavity. Since we chose to develop these results for graph
motion with a free boundary, although the general lines of proof
(via maximum principles) have precedents, many details had to be
developed anew. For example, section 13 contains an extension of the
maximum principle for symmetric tensors with Neumann-type boundary
conditions given in \cite{Stahl}, which in our setting allows one to
show preservation of weak concavity in general. The results obtained
in sections 11-14 are summarized in the following theorem.\newpage

\textbf{Theorem 1.2.} Let $T_{max}$ be the maximal existence time
for the evolution. Assuming $T_{max}<\infty$, the second fundamental
form $h$ is unbounded at the junction $\Gamma_t$, as $t\rightarrow
T_{max}$:
$$\limsup_{t\rightarrow T_{max}}(\sup_{\Gamma_t}|h|_g)=\infty.$$
If the mean curvature of the initial hypersurface is strictly
negative ($\sup_{\Sigma_0}H=H_0<0$), then $T_{max}$ is finite. If
$\Sigma_0$ is weakly concave ($h\leq 0$ at $t=0$), this is preserved
by the evolution. \vspace{.2cm}

The expected global existence result is that, assuming weak
concavity, $diam(\Sigma_t)\rightarrow 0$ as $t\rightarrow
T_{max}.$\vspace{.2cm}

\emph{Acknowledgments.} It is a pleasure to thank Nicholas Alikakos
for originally proposing to consider the problem of mean curvature
flow for networks of surfaces meeting at constant angles, and for
his interest in this work. Most of the work on short-time existence
was undertaken during a stay at the Max-Planck Institute for
Gravitational Physics in Golm (January-June, 2007); I am grateful to
the Max-Planck Society for supporting the visit, and to Gerhard
Huisken, director of the Geometric Analysis group, for the
invitation. Finally, thanks to Mariel S\'{a}ez for communicating the
results of the Lens Seminar (\cite{LensSeminar}) and of her recent
work on mean curvature flow of networks (partly in collaboration
with Rafe Mazzeo, \cite{MazzeoSaez}).\vspace{.4cm}

\textbf{2. Normal velocity of the moving boundary.} The evolution is
naturally supplied with initial data $\Sigma_0$, a graph meeting
$\mathbb{R}^{n+1}$ at the prescribed angle. Since we are interested
in classical solutions in the parabolic H\"{o}lder space
$C^{2+\alpha,1+\alpha/2}$, we expect an additional compatibility
condition at $t=0$. We discuss this first for graph m.c.m. $w(y,t)$.

Denote by $\Gamma(t)$ a global parametrization of $\partial D(t)$
(with domain in a fixed manifold, and `space variables' left
implicit). Differentiating in $t$ the `contact condition'
$w(\Gamma(t),t)=0$, we find:
$$w_t+\langle Dw,\dot{\Gamma}(t)\rangle =0.$$
Denote by $n_t$ the unit normal vector field to $\Gamma(t)$, chosen
so that $\langle n_t,Dw\rangle>0$. The contact condition also
implies the gradient of $w$ is purely normal:
$$Dw_{|\partial D(t)}=(D_{n_t}w)n_t.$$
Combining this with the angle condition, and bearing in mind that
$D_{n_t}w_{|\Gamma(t)}>0$, we find:
$$D_{n_t}w=\frac {\beta_0}{\beta}\mbox{ on }\partial D(t),\quad \beta_0:=\sqrt{1-\beta^2}.$$
Thus, on $\partial D(t)$:
$$\frac 1{\beta}H=\sqrt{1+(D_{n_t}w)^2}H=w_t=-\langle
\dot{\Gamma}(t),n_t\rangle D_{n_t}w=-\dot{\Gamma}_n(t)\frac
{\beta_0}{\beta},$$ and we find the normal velocity of the moving
boundary (independent of the parametrization of $\Gamma_t$):
$$\dot{\Gamma}_n=-\frac 1{\beta_0}H_{|\Gamma(t)},$$
which in particular must hold at $t=0$. Note that we don't get a
`compatibility condition' in the usual sense (of a constraint on the
2-jet of the initial data), but instead an equation of motion for
the moving boundary. (Later, in the fixed-domain formulation, we
will have to deal with a real compatibility condition).
\vspace{.2cm}

Now consider mean curvature flow parametrized over a time-dependent
domain ${\cal D}(t)$, with the boundary conditions:
$$\langle F(x,t),e_{n+1}\rangle=0,\quad \langle N,e_{n+1}\rangle=\beta,\quad
x\in \partial{\cal D}(t).$$
 Let $\nu=\nu_t$
be the inner unit normal to ${\cal D}(t)$. Suppose ${\cal
S}(\theta,t)$, $\theta \in S^{n-1}$, parametrizes $\partial {\cal
D}(t)$; thus the `junction' $\partial \Sigma_t$ is parametrized by
$\Gamma(\theta,t)=F(t,{\cal S}(\theta,t))$, and (denoting partial
$t$ derivatives with a dot):
$$\dot{\Gamma}(\theta,t)=\partial_tF(t,{\cal S}(\theta,t))+dF[{\cal
S}(\theta,t)]=HN+(\dot{\cal S}\cdot \nu)\partial_{\nu}F$$ (where we
used the fact that $\partial_{\tau}F=0$ for any $\tau\in
T\partial{\cal D}(t)$.) Thus, using $\langle N,n\rangle=-\beta_0$:
$$\dot{\Gamma}_n:=\dot{\Gamma}(\theta,t)\cdot
n=-\beta_0H+(\dot{\cal S}\cdot \nu)\langle
\partial_{\nu}F,n\rangle.$$
On the other hand, from $\langle F(t,{\cal
S}(\theta,t)),e_{n+1}\rangle\equiv 0$, we find by differentiation:
$$H\beta+(\dot{\cal S}\cdot \nu)\langle
\partial_{\nu}F,e_{n+1}\rangle=0,$$
or $\dot{\cal S}\cdot \nu=-H\beta/\langle
\partial_{\nu}F,e_{n+1}\rangle$. Letting
$T:=\frac{\partial_{\nu}F}{|\partial_{\nu}F|}$ (tangent to
$\Sigma_t$ at the interface), we have:
$$\dot{\Gamma}_n=-H(\beta_0+\beta\frac{\langle T,n\rangle}{\langle
T,e_{n+1}\rangle}).$$ Denoting by $N'=N-\beta e_{n+1}$ the
$\mathbb{R}^n$ component of $N$, we clearly have $n=-(1/\beta_0)N'$,
so $\langle T,n\rangle=-(1/\beta_0)\langle
T,N'\rangle=(\beta/\beta_0)\langle T,e_{n+1}\rangle$, and we
conclude:
$$\dot{\Gamma}_n=-H(\beta_0+\frac{\beta^2}{\beta_0})=-\frac
1{\beta_0}H,$$ as before.\vspace{.2cm}

\emph{Remark 2.1.} This is not unexpected, if we accept there is a
reparametrization connecting the two motions, respecting the
boundary conditions. That is, the ODE argument in
\cite{EckerHuisken} should also work in the presence of boundary
conditions and moving boundaries.\vspace{.2cm}

\emph{Remark 2.2.} We remark that for more general (non-symmetric,
non-flat) triple junctions with 120 degree angles, the condition:
$$H^1+H^2=H^3\mbox{ on }\Gamma(t)$$
must hold at the junction (for graphs, oriented by the upward
normal), which in particular gives a geometric constraint on the
initial data, for classical evolution in $C^{2+\alpha, 1+\alpha/2}$.
This is \emph{automatic} in the symmetric case ($w^2=-w^1$), since
$H^3=0$ and $H^I=tr_{g^I}d^2w^I$ for $I=1,2$. \vspace{.3cm}

\textbf{3. Choice of `gauge'.} It is traditional in moving boundary
problems to parametrize the time-dependent domain $D(t)$ of the
unknown $w(y,t)$ by a time-dependent diffeomorphism:
$$y=\varphi (x,t),\quad \varphi(t):D_0\rightarrow D(t),$$
and then derive the equation satisfied by the coordinate-changed
function from the equation for $w$ (see e.g. \cite{LunardiBaconneau}
or \cite{Solonnikov}). Motivated by the work on curve networks
(\cite{Mantegazza et al.}) we will, instead, consider a general
parametrization:
$$F:D_0\times [0,T]\rightarrow {\mathbb R}^{n+1},\quad
F(x,t)=[\varphi(x,t),u(x,t)]\in\mathbb{R}^n\times \mathbb{R}$$ and
derive an equation for $F$ directly from the definition of mean
curvature motion:
$$\langle \partial_tF,N\rangle=H.$$
(We'll still assume $\varphi(t):D_0\rightarrow D(t)$ is a
diffeomorphism.) The first and second fundamental forms are given
by:
$$g_{ij}=\langle F_i,F_j\rangle, \quad A(F_i,F_j)=\langle
F_{ij},N\rangle.$$ (Notation: $DF=F_ie_i, D^2F(e_i,e_j)=F_{ij}$,
$(e_i)$ is the standard basis of $\mathbb R^{n+1}$.) The mean
curvature is the trace of $A$ in the induced metric:
$$H=\langle g^{ij}(DF)F_{ij},N\rangle.$$ The equation for $F$ is:
$$\langle \partial_tF-g^{ij}(DF)F_{ij},N\rangle=0.$$

There is a natural `gauge choice' yielding a quasilinear parabolic
system:
$$\partial_tF-g^{ij}(DF)F_{ij}=0.$$
We will sometimes refer to this as the `split gauge', since in terms
of the components $F=[\varphi,u]$ we have the essentially decoupled
system:
$$\left \{ \begin{array}{ccc}
\partial_tu-g^{ij}(D\varphi,Du)u_{ij}&=&0,\\
\partial_t\varphi-g^{ij}(D\varphi,Du)\varphi_{ij}&=&0.\end{array} \right .$$ The
splitting is useful to state the boundary conditions:
\begin{eqnarray*}
&u_{|\partial D_0}=0\mbox{ (`contact')},\\
&N^{n+1}(D\varphi,Du)_{|\partial D_0}=\beta\mbox{ (`angle')}.
\end{eqnarray*}
 We immediately see there is a problem,
since we have 2 scalar boundary conditions for $n+1$ unknowns (and
no moving boundary to help!) Our solution to this is to introduce
$n-1$ additional `orthogonality conditions' at the boundary for the
parametrization $\varphi(t)$. We impose:
$$\langle D_{\tau}\varphi,D_n\varphi\rangle_{|\partial D_0}=0,$$
for any $\tau\in T\partial D_0$, where $n$ is the inward unit normal
to $D_0$.

Geometrically, the `orthogonality' boundary condition has precedent
in a method often adopted when dealing with the evolution of
hypersurfaces in $\mathbb{R}^{n+1}$ intersecting a fixed
$n$-dimensional `support surface' orthogonally (see e.g.
\cite{Struwe}): one replaces vanishing inner product of the unit
normals (a single scalar condition) by a stronger Neumann-type
condition for the parametrization, corresponding to $n-1$ scalar
conditions. (More details are given in Section 10.)\vspace{.2cm}

The system must also be supplied with initial data. We assume given
an initial hypersurface $\Sigma_0$, the graph of a $C^{3+\alpha}$
function $\tilde u_0(x)$ defined in the $C^{3+\alpha}$ domain
$D_0\subset \mathbb{R}^n$. (The reason for this choice of
differentiability class will be seen later.) It would seem natural
to set $\varphi_0=Id_{D_0}$, but this causes problems (related to
compatibility; see Section 4 below). We do require the 1-jet of
$\varphi_0$ at the boundary to be that of the identity:
$${\varphi_0}_{|\partial D_0}=Id,\quad D{\varphi_0}_{|\partial D}
=\mathbb{I}.$$ (In particular, the orthogonality condition holds at
$t=0$.)

We need a more explicit expression for the unit normal, and for that
we use the `vector product':
$$\tilde{N}(D\varphi,Du):=(-1)^n\det\left
[\begin{array}{ccc}e_1&\cdots&e_{n+1}\\DF^1&\cdots&DF^{n+1}\end{array}
\right]=(-1)^{n}\det
\left[\begin{array}{cccc}e_1&\ldots&e_n&e_{n+1}\\D\varphi^1&\ldots&D\varphi^n&Du\end{array}\right]$$
$$:=[J(D\varphi,Du),J_{\varphi}]\in {\mathbb R}^n\times
{\mathbb R},$$
 where $DF^i\in {\mathbb R}^n$ for $i=1,\ldots n+1$, $J_{\varphi}>0$ is the jacobian of $\varphi$ and
$(-1)^n$ is introduced to make sure the last component is positive.
$J(D\varphi,Du)$ is an ${\mathbb R}^n$-valued multilinear form,
linear in the components $u_i$ of $Du$ and of weight $n-1$ in the
components of $D\varphi$. It is easy to check that
$J(\mathbb{I},Du)=-Du$. The unit normal is:
$$N(D\varphi,
Du)=\tilde{N}(D\varphi,Du)/(|J(D\varphi,Du)|^2+(J_{\varphi})^2)^{1/2}.$$
Thus the angle condition may be stated in the form:
$$\beta[|J(Du,D\varphi)|^2+(J_{\varphi})^2]^{1/2}_{|\partial
D_0}={J_{\varphi}}_{|\partial D_0},$$ and we lose nothing by
squaring it:
$$B(D\varphi,Du):=
\beta^2|J(Du,D\varphi)|^2-\beta_0^2(J_{\varphi})^2_{|\partial
D_0}=0.$$ \vspace{.3cm}

\textbf{4. Compatibility and the choice of $\varphi_0$.} Assume
$D{\varphi_0}_{|\partial D_0}=\mathbb{I}$. Differentiating in $t$
the contact condition $u_{|\partial D_0}=0$ and evaluating at $t=0$,
we find:
$$0=g^{ij}(\mathbb I,Du_0)u_{0ij}\equiv g_0^{ij}u_{0ij}\mbox{ on
}\partial D_0.$$ To interpret this condition, consider the mean
curvature at $t=0$, on $\partial D_0$:
$$H_0=\frac 1{v_0}[\langle
J(\mathbb{I},Du_0),g_0^{ij}\varphi_{0ij}\rangle+J_{\varphi_0}g^{ij}u_{0ij}],$$
where:
$$v_0=[|J(\mathbb{I},Du_0)|^2+J_{\varphi_0}^2]^{1/2}_{|\partial D_0}
=(|Du_0|^2+1)^{1/2}_{|\partial D_0}=\frac 1{\beta},$$ using (recall
$\beta_0:=\sqrt{1-\beta^2}$):
$$J(\mathbb I,Du_0)=-Du_0=-(D_nu_0)n=\frac
{\beta_0}{\beta}n$$ on $\partial D_0$. Thus the compatibility
condition is equivalent to:
$${H_0}_{|\partial D_0}=-\beta_0 g_0^{ij}\langle
\varphi_{0ij},n\rangle_{|\partial D_0}.$$ This implies we can't
choose $\varphi_0\equiv Id$ (on all of $D_0$), unless
${H_0}_{|\partial D_0}\equiv 0$, a constraint not present in the
geometric problem (as seen above). Instead, regarding $H_0$ as given
(by $\Sigma_0$), and using:
$$g_0^{ij}=\delta_{ij}-\frac{u_{0i}u_{0j}}{v_0^2}=\delta_{ij}-\beta_0^2n^in^j,$$
we find the compatibility constraint:
$$\langle
(\delta_{ij}-\beta_0^2n^in^j)\varphi_{0ij},n\rangle=-\frac
1{\beta_0}H_0\mbox{ on }\partial D_0.$$ Given the zero and first
order constraints on $\varphi_0$, this can also be written as:
$$n^in^j\langle \varphi_{0ij},n\rangle=-\frac
1{\beta^2\beta_0}H_0\mbox{ on }\partial D_0.$$ The next lemma shows
this can be solved.\vspace{.2cm}

\textbf{Lemma 4.1.} Let $D_0\subset \mathbb{R}^n$ be a uniformly
$C^{3+\alpha}$ domain (possibly unbounded), $h\in
C^{\alpha}(\partial D_0)$ $(0<\alpha<1)$.

(i) One may find a diffeomorphism $\varphi\in Diff^{2+\alpha}(D_0)$
satisfying on $\partial D_0$:
$$\varphi=Id, \quad d\varphi=\mathbb{I},\quad n\cdot
d^2\varphi(n,n)=h.$$

(ii) More generally, given a non-vanishing vector field $e\in
C^{1+\alpha}(\partial D_0;{\mathbb R}^n)$, one may find $\varphi\in
Diff^{2+\alpha}(D_0)$ satisfying on $\partial D_0$:
$$\varphi=Id, \quad d_n\varphi=e, \quad n\cdot
d^2\varphi(n,n)=h.$$ If $\partial D_0$ has two components, we may
even require $\varphi$ to satisfy the conditions in parts (i) and
(ii) at $\partial_1D_0$, $\partial_2D_0$ (resp.), with different
functions $h$. (This will be needed in section 11).
 \vspace{.1cm}

As usual, a domain is `uniformly $C^{3+\alpha}$' if at each boundary
point there are local charts to the upper half-space (of class
$C^{3+\alpha}$), defined on balls of uniform radius, and with
uniform bounds on the $C^{3+\alpha}$ norms of the charts and their
inverses.\vspace{.2cm}

\emph{Remarks:} 4.1.The proof is given in Appendix 1.

4.2. Note that, in particular, $\varphi$ satisfies the orthogonality
conditions at $\partial D_0$.

4.3. It is at this step in the proof that we have a drop in
regularity: for $C^{2+\alpha}$ local solutions, we require
$C^{3+\alpha}$ initial data. While this is not unexpected in
free-boundary problems (see e.g. \cite{LunardiBaconneau}), I don't
know a counterexample to the lemma if $D_0$ is assumed to be a
$C^{2+\alpha}$ domain.

4.4. In our application of the lemma, we in fact have $h\in
C^{1+\alpha}(\partial D_0)$, but this does not imply higher
regularity for $\varphi$. \vspace{.3cm}

\textbf{5. Linearization.} The evolution equation and boundary
conditions in `split gauge' are: $$\left \{\begin{array}{ccc}
 F_t-g^{ij}(DF)F_{ij}&=&0,\\
 u_{|\partial D_0}&=&0,\\
 B(D\varphi,Du)_{|\partial
D_0}&=&0,\\
{\cal O}(D\varphi)_{|\partial D_0}&=&0,
\end{array}\right .$$
where: $$ {\cal O}(D\varphi):=\langle
D^T\varphi,D_n\varphi\rangle.$$ Here
$D^T\varphi=D\varphi-(D_n\varphi)\langle \cdot,n\rangle$ is an
${\mathbb R}^n$-valued 1-form on $\partial D_0$. We'll prove
short-time existence for this system (with initial data
$u_0,\varphi_0$) in $C^{2+\alpha, 1+\alpha/2}$ by the usual
fixed-point argument based on linear parabolic theory. Given
$\bar{F}=[\bar{\varphi},\bar{u}]$ in a suitable ball in this
H\"older space with center $F_0=[\varphi_0,u_0]$, it suffices to
consider the `pseudolinearization' of the system:
$$F_t-g^{ij}(DF_0)u_{ij}=[g^{ij}(D\bar{F})-g^{ij}(DF_0)]\bar{F}_{ij}:={\cal
F}(\bar{F},\bar{F_0});$$ a fixed point of the map $\bar{F}\mapsto F$
corresponds to a solution of the quasilinear equation.\vspace{.2cm}

For the nonlinear boundary conditions, we need the honest
linearization at $F_0$. For the angle condition, a computation using
the boundary constraints on $u_0$ and $\varphi_0$ yields:
$$\frac 12 {\cal
L}_0B[D\varphi,Du]=-\beta\beta_0D_nu-\beta_0^2\langle
D_n\varphi,n\rangle.$$ The corresponding linear boundary condition
will be:
$$\beta\beta_0D_nu+(1-\beta^2)\langle
D_n\varphi,n\rangle={\cal B}(D\bar F,DF_0),$$ where:
$$2{\cal B}(D\bar F,DF_0):=B(D\bar{\varphi},D\bar
u)-B(Du_0,D\varphi_0)-{\cal
L}_0B[D(\bar{\varphi}-\varphi_0),D(\bar{u}-u_0)],$$ and we used:
$$-\frac 12 {\cal L}_0[D\varphi_0,Du_0]_{|\partial
D_0}=\beta\beta_0D_nu_0+(1-\beta^2)\langle
D_n\varphi_0,n\rangle_{|\partial D_0}=0.$$ Also,
$B(D\varphi_0,Du_0)_{|\partial D_0}=0$, so at a fixed point
$B(D\varphi,Du)_{|\partial D_0}=0$.\vspace{.2cm}

Linearizing the orthogonality boundary condition, we find that
${\cal L}_0{\cal O}[D\varphi]$ is the 1-form on $\partial D_0$:
$${\cal
L}_0{\cal
O}[D\varphi](v)=(\partial_j\varphi^i+\partial_i\varphi^j)n^j(\delta_{ik}-n^kn^i)v^k$$
(with sum over repeated indices.) The corresponding linear boundary
condition is:
$$\langle D_n\varphi,proj^T\rangle+\langle
D^T\varphi,n\rangle=-\Omega(D\bar{\varphi},D\varphi_0),$$ where:
$$\Omega(D\bar{\varphi},D\varphi_0):={\cal O}(D\bar{\varphi})-{\cal O}(D\varphi_0)-{\cal
L}_0{\cal O}[D\bar{\varphi}-D\varphi_0],$$ and we used:
$${\cal L}_0{\cal O}[D\varphi_0]_{|\partial D_0}=\langle
(D_n\varphi_0)^T,\cdot\rangle+\langle
D^T\varphi_0,n\rangle_{|\partial D_0}=0.$$ \vspace{.3cm}

\textbf{6. Complementarity.} We wish to apply linear existence
theory to the system:
$$F_t-g^{ij}(DF_0)F_{ij}=\bar{\cal F},$$
with boundary conditions at $\partial D_0$:
$$u=0,$$
$$\beta\beta_0D_nu+\beta_0^2\langle D_n\varphi,n\rangle
=\bar{\cal B},$$
$$\langle D_n\varphi,proj^T\rangle +\langle
D^T\varphi,n\rangle=-\bar{\Omega}$$ and initial conditions:
$$u_{t=0}=u_0,\quad \varphi_{t=0}=\varphi_0.$$
It is easy to see that the initial data satisfy the linearized
boundary conditions, and above we constructed $\varphi_0$ so as to
guarantee $g^{ij}(Du_0,D\varphi_0)u_{0ij}{_{|\partial D_0}}=0.$
(There is no first-order compatibility condition for $\varphi_0$.)
Thus the linear system satisfies the required compatibility at
$t=0$.\vspace{.2cm}

Since the linearized boundary conditions are slightly non-standard,
we must verify they satisfy the `complementarity'
(Lopatinski-Shapiro) conditions. We \emph{fix} $x_0\in \partial D_0$
and introduce adapted coordinates $(\rho,\sigma)$ in a neighborhood
${\cal N}_0\subset {\cal N}$ of $x_0$ in $D_0$:
$$x=\Gamma_0+\rho n(\sigma),\quad \sigma=(\sigma_a)\in {\cal U},$$
where $\Gamma_0:{\cal U}\rightarrow {\mathbb R^n}$ is a local chart
for $\partial D_0$ at $x_0$ (${\cal U}\subset {\mathbb R}^{n-1}$
open). This defines a basis of tangential vector fields in
$\Gamma_0({\cal U})$, and we may assume that, at $x_0$: $\langle
\tau_a,\tau_b\rangle=\delta_{ab}$ and
$\nabla_{\tau_a}\tau_b(x_0)=0$. Let $U$ and $\psi$ be defined in
$(-\rho_1,0)\times {\cal U}\times [0,T]$ by:
$$U(\rho,\sigma,t)=u(\Gamma_0(\sigma)+\rho n(\sigma),t),\quad
\psi(\rho,\sigma,t)=\varphi(\Gamma_0(\sigma)+\rho n(\sigma),t).$$ In
these coordinates, the induced metric is written (in `block form'):
$$[g]=\left [ \begin{array}{cc}|\psi_{\rho}|^2+(U_{\rho})^2 &\langle
\psi_{\rho},\psi_{a}\rangle+U_{\rho}U_a\\\langle
\psi_{\rho},\psi_a\rangle+U_{\rho}U_a &\langle
\psi_{a},\psi_b\rangle+U_aU_b\end{array}\right ]=\left
[\begin{array}{cc}\frac 1{\beta^2}& 0\\0& \mathbb{I}_{n-1}
\end{array}\right ]$$ at $t=0$ and $x_0$.

We have: $$U_{\rho \rho}=D^2u(n,n)\mbox{ (since $\nabla_nn=0$)},$$

$$U_{ab}=D^2u(\tau_a,\tau_b)+Du\cdot
\nabla_{\tau_a}\tau_b=D^2u(\tau_a,\tau_b) \mbox{ at }x_0,$$ and we
don't need $U_{\rho a}$, since $g_{\rho a}=0$ at $x_0$.

Thus: $$tr_{g_0}D^2u(x_0)=\beta^2 D^2u(n,n)+\sum_a
D^2u(\tau_a,\tau_a)=\beta^2U_{\rho \rho}+\sum_aU_{aa}:=
\beta^2U_{\rho \rho}+\Delta_{\sigma}U,$$ and, likewise:
$$tr_{g_0}D^2\varphi(x_0)=\beta^2\psi_{\rho
\rho}+\Delta_{\sigma}\psi.$$ For the linearized orthogonality
operator, note that, at $x_0$:
$${\cal L}_0{\cal O}[D\psi]=(\langle
\psi_{\rho},\tau_a\rangle+\psi_a,n\rangle )\tau_a.$$ Putting
everything together, the linear system to consider at $x_0$ is:
$$U_t-\beta^2U_{\rho \rho}-\Delta_{\sigma}U=0,$$
$$\psi_t-\beta^2\psi_{\rho \rho}-\Delta_{\sigma}\psi=0,$$
with boundary conditions: $U|_{\rho=0}=0$,
$$\beta_0\langle \psi_{\rho},n\rangle+\beta
U_{\rho}|_{\rho=0}=b(\sigma,t),$$
$$\langle \psi_{\rho},\tau_a\rangle+\langle
\psi_a,n\rangle|_{\rho=0}=\omega_a(\sigma,t),\quad a=1,\ldots n-1.$$
\vspace{.2cm}

Now take Fourier transform in $\sigma\in {\mathbb R}^{n-1}$, Laplace
transform in $t$ to obtain:
$$\hat{U}(\rho,\xi, p)\in \mathbb{C}, \hat{\psi}(\rho,\xi,p)\in
\mathbb{C}^n; \quad \xi\in{\mathbb R}^{n-1},p\in
\mathbb{C},\rho<0.$$ In transformed variables, we obtain the system
of linear ODE (in $\rho<0$, for fixed $(\xi,p)$):
$$\beta^2\hat{U}_{\rho \rho}-(p+|\xi|^2)\hat{U}=0,$$
$$\beta^2\hat{\psi}_{\rho \rho}-(p+|\xi|^2)\hat{\psi}=0.$$
Writing the solution in the form:
 $$\left[\begin{array}{c}\hat{U(\rho)}\\\hat{\psi(\rho)}\end{array}\right]
 =e^{i\rho\gamma}\left[\begin{array}{c}\hat{U}(0)\\\hat{\psi}(0)\end{array}\right],$$
 we find the characteristic equation $\beta^2\gamma^2+p+|\xi|^2=0$,
 and choose the root $\gamma$ so that
 $i\gamma=(1/\beta)\sqrt{\Delta}$ (where $\Delta=p+|\xi|^2$ and we
 take the branch of $\sqrt{}$: $Re(\sqrt{\Delta})>0$). Here $(p,\xi)\in {\cal A}$, where:
 $${\cal A}=\{(p,\xi)\in {\mathbb C}\times \mathbb {R}^{n-1};
 |p|+|\xi|>0, Re(p)>-|\xi|^2\}.$$ Thus the solutions decay as
 $\rho\rightarrow -\infty$. Let ${\cal W}^+$ be the space of such
 decaying solutions, dim$_{\mathbb C}{\cal W}^+=n-1$. The relevant
 boundary operator on ${\cal W}^+$ is:
 $${\mathbb B}\left[\begin{array}{c}\hat{U}\\\hat{\psi}\end{array}\right]
 =\left[\begin{array}{c}\hat{U}\\ \beta_0\langle
 \hat{\psi}_{\rho},n\rangle+\beta \hat{U}_{\rho}\\
 \langle \hat{\psi}_{\rho},\tau_a\rangle+i\xi_a\langle
 \hat{\Psi},n\rangle\end{array}\right]_{|\rho=0}=
 \left[\begin{array}{c}\hat{U}(0)\\\beta_0(i\gamma)\langle
 \hat{\psi}(0),n\rangle +i\beta\gamma\hat{U}(0)\\
 (i\gamma)\langle \hat{\psi}(0),\tau_a\rangle+i\xi_a\langle
 \hat{\psi}(0),n\rangle \end{array}\right]$$(a vector in ${\mathbb
 C}\times{\mathbb C}\times{\mathbb C}^{n-1}).$

 The `complementarity condition' (see e.g. \cite{EidelmanZhitarasu}) is the statement that $\mathbb B$
 is a linear isomorphism from ${\cal W}^+$ to $\mathbb{C}^{n+1}$.
 With respect to the basis $\{\hat{U}(0),\langle
 \hat{\psi}(0),n\rangle,\langle \hat{\psi}(0)\tau_a\rangle\}$ of ${\cal
 W}^+$, the matrix of $\mathbb B$ is (in `block form'):
 $$[\mathbb{B}]=\left [\begin{array}{ccc}1&0&[0]_{1\times
 (n-1)}\\\sqrt{\Delta}&\frac{\beta_0}{\beta}\sqrt{\Delta}&[0]_{1\times
 (n-1)}\\ \left [0\right ]_{(n-1)\times 1}&[i\xi_a]_{(n-1)\times
 1}&\frac{\sqrt{\Delta}}{\beta}\mathbb{I}_{n-1}\end{array}\right].$$
 This is triangular with non-zero diagonal entries for every
 $(p,\xi)\in {\cal A}$. Hence $\mathbb{B}$ is an isomorphism.

 \vspace{.4cm}
 \textbf{7. Estimates in H\"{o}lder spaces.}

 For the fixed-point argument based on the linear system, we need
 estimates for $||\cal F||_{\alpha}$, $||{\cal B}||_{1+\alpha}$,
 $||\Omega||_{1+\alpha}$, of two types: `mapping' and `contraction'
 estimates.\vspace{.2cm}

 A bit more precisely, for $T>0$, $R>0$ and $Q^T=D_0\times [0,T]$ consider the open ball:
 $$B_R^T=\{F\in C^{2+\alpha,1+\alpha/2}(Q^T,{\mathbb
R}^{n+1});||F-F_0||_{2+\alpha}<R, F|_{t=0}=F_0\}.$$
($F_0=[\varphi_0,u_0]$ is defined from the initial surface
$\Sigma_0$, via Lemma 4.1.) Solving the linear system with
`right-hand side' defined by $\bar{F}\in B_R^T$ defines a map
${\mathbb F}: \bar{F}\mapsto F$, and we need to verify that, for
suitable choices of $T$ and $R$, $\mathbb{F}$ maps into $B_R^T$ and
is a contraction.\vspace{.2cm}

\emph{Remark:} The argument that follows is standard, and the
experienced reader may want to skip to the statement of local
existence at the end of the next section. On the other hand the
result is not covered by any general theorem proved in detail in a
reference known to the author, and some readers may find it useful
to have all the details included. Another reason is that, although
the `right hand sides' are clearly quadratic, without explicit
expressions one might run into trouble with compositions (which
behave poorly in H\"{o}lder spaces), or when appealing to `Taylor
remainder arguments' if the domain is not convex.\vspace{.2cm}

 For `mapping', we need estimates of the form:
 $$||{\cal F}(\bar{F},F_0)||_{\alpha}+||{\cal
 B}(D\bar{F},DF_0)||_{1+\alpha}+||\Omega(D\bar{\varphi},D\varphi_0)||_{1+\alpha}\mbox{
 decays as }T\rightarrow 0_+,$$
 and for `contraction':
 $$||{\cal F}(F^1,F^0)-{\cal F}(F^2,F^0)||_{\alpha}+||{\cal
 B}(DF^1,DF^2)||_{1+\alpha}+||\Omega(D\varphi^1,D\varphi^2)||_{1+\alpha}\leq
 \mu(T)||F^1-F^2||_{2+\alpha},$$
 where $\mu(T)\rightarrow 0$ as $T\rightarrow 0_+$.\vspace{.2cm}

 \emph{Notation:} The
 $(\alpha,\alpha/2)$ norms are taken on $Q^T$, the
 $(1+\alpha,(1+\alpha)/2)$ norms on $\partial D_0\times [0,T]$).
 Double bars without an index refer to the $(2+\alpha,1+\alpha/2)$
 norm, single bars to supremum norms over $Q^T$, and parabolic norms are
 indexed by their spatial regularity ($\alpha$ for
 ($\alpha,\alpha/2$), etc.) In general, we use brackets for
 H\"{o}lder-type difference quotients.
  \vspace{.2cm}

 We deal with the estimates for the `forcing term' $\cal F$ first. Consider the map $${\cal
 G}:Imm(\mathbb{R}^n,\mathbb{R}^{n+1})\rightarrow GL_n$$ which
 associates to the linear immersion $A$ the inverse matrix of
 $(\langle A_i,A_j\rangle)_{i=1}^n$, inner products of the rows of
 $A$. ${\cal G}$ is smooth, in particular locally Lipschitz in the
 space ${\cal W}$ of linear immersions. Hence if $F^1,F^2$ are maps
 $Q^T\rightarrow \mathbb{R}^{n+1}$, such that
 $DF^i\in C^{\alpha, \alpha/2}(Q^T)$ and $DF^i(z)\in K$ for all
 $z\in Q^T$, where $K\subset {\cal W}$ is a fixed compact set, we
 have the bound:
 $$||{\cal G}(DF^1)-{\cal G}(DF^2)||_{\alpha}\leq
 c_K||D(F^1-F^2)||_{\alpha}.$$
 In fact our maps $F^i$ are in $C^{2+\alpha,1+\alpha/2}$, so
 $DF^i\in C^{1+\alpha,\frac{1+\alpha}2}$. From this higher
 regularity we obtain the decay as $T\rightarrow 0_+$.
 Assuming $F^1|_{t=0}=F^2|_{t=0}$, we have:
 $$|D(F^1-F^2)|\leq
 [D(F^1-F^2)]_t^{(\frac{1+\alpha}2)}T^{\frac{1+\alpha}2}.$$

 To continue, we recall an elementary fact for H\"{o}lder
 spaces:\vspace{.2cm}

  Let $D\subset {\mathbb R}^n$ be a uniformly
 $C^1$ domain (not necessarily convex or bounded). Then if $f\in
 C^1(D)$ and $\alpha\in (0,1)$, we have: $$[f]^{(\alpha)}\leq
 C_D||f||_{C^1}.$$\vspace{.2cm}

 Here `uniformly $C^1$' means $D$ can be covered by countably many
 balls of a fixed radius, which are domains of $C^1$
 manifold-with-boundary local charts for $D$, with uniform $C^1$
 bounds for the charts and their inverses. The constant $C_D$
 depends on those bounds. Applying the lemma to $DF$, where
 $F=F^1-F^2$ vanishes identically at $t=0$, and assuming $T<1$:

 $$[DF]_x^{(a)}\leq c(|DF|+|D^2F|)\leq
 c([DF]_t^{(\frac{1+\alpha}2)}T^{\frac{1+\alpha}2}+[D^2F]_t^{(\frac{\alpha}2)}T^{\alpha/2})\leq
 c||F||T^{\alpha/2},$$
 (where $c$ depends on $D_0$) and similarly for the oscillation in $t$:
 $$[DF]_t^{(\frac{\alpha}2)}\leq
 [DF]_t^{(\frac{1+\alpha}2)}T^{1/2}\leq ||F||T^{1/2},$$
 so we have:
 $$||D(F^1-F^2)||_{\alpha}\leq c ||F^1-F^2||T^{\alpha/2}.$$
 We conclude, under the
 assumption $F^1=F^2$ at $t=0$:
 $$||{\cal G}(DF^1)-{\cal G}(DF^2)||_{\alpha}\leq
 c_K||F^1-F^2||T^{\alpha/2}.$$
 In particular, applying this to $\bar{F}$ and $F_0$, we find:
 $$||({\cal G}(D\bar{F})-{\cal G}(DF_0))\bar{D^2F}||_{\alpha}\leq
 c_K||\bar{F}-F_0||T^{\alpha/2}||\bar{F}||,$$
 and for $F^1$ and $F^2$ coinciding at $t=0$:
 $$||({\cal G}(DF^1)-{\cal G}(DF^2))D^2F^1||_{\alpha}\leq
 c_K||F^1-F^2||T^{\alpha/2}||F^1||,$$
 as well as:
 $$||({\cal G}(DF^2)-{\cal G}(DF_0))(D^2F^1-D^2F^2)||_{\alpha}\leq
 c_K||F^2-F_0||T^{\alpha/2}||F^1-F^2||,$$
 so we have the mapping and contraction estimates for ${\cal
 F}(\bar{F},F_0)$ and ${\cal F}(F^1,F_0)-{\cal
 F}(F^2,F_0)$.

 \textbf{Lemma 7.1.} Assume $\bar{F},F_0,F^1,F^2$ are in
 $C^{2+\alpha,1+\alpha/2}(Q^T;{\mathbb R}^{n+1})$ and have the same
 initial values, and that
 $D\bar{F},DF_0,DF^1,DF^2$ all take
 values in the compact subset $K$ of $Imm(\mathbb{R}^n,\mathbb
 {R}^{n+1})$. Then:
 $$||{\cal F}(\bar{F},F_0)||_{\alpha}\leq
 c_K||\bar{F}-F_0||||\bar{F}||T^{\alpha/2},$$
 $$||{\cal F}(F^1,F_0)-{\cal F}(F^2,F_0)||_{\alpha}\leq
 c_K(||F^1||+||F^2-F_0||)T^{\alpha/2}||F^1-F^2||.$$ In particular,
 if $\bar{F}\in B_R^T$:
 $$||{\cal F}(\bar{F},F_0)||_{\alpha}\leq c_0RT^{\alpha/2}.$$
 If $\bar{F}^1,\bar{F}^2\in B_R^T$, we have:
$$||{\cal F}(\bar{F}^1,F_0)-{\cal F}(\bar{F}^2,F_0)||_{\alpha}\leq
c_0T^{\alpha/2}||\bar{F}^1-\bar{F}^2||.$$
 (The constant $c_0$ depends only on the data at $t=0$, and we assume
 $T<1$, $R<1$).\vspace{.2cm}

 Turning to the orthogonality boundary condition, first observe
 that:
 $$\Omega(D\varphi^1,D\varphi^2)=\langle
 D^T\varphi^1,D_n\varphi^1\rangle-\langle
 D^T\varphi^2,D_n\varphi^2\rangle-{\cal L}_0{\cal
 O}[D\varphi^1-D\varphi^2]$$
 $$=\langle
 D^T(\varphi^1-\varphi^2),D_n\varphi^1\rangle+\langle D^T\varphi^2,D_n(\varphi^1-\varphi^2)\rangle-\langle
 D_n(\varphi^1-\varphi^2),D^T\varphi_0\rangle-\langle
 D^T(\varphi^1-\varphi^2),D_n\varphi_0\rangle$$
 $$=\langle
 D^T\varphi^1-D^T\varphi^2,D_n\varphi^1-D_n\varphi_0\rangle
 +\langle
 D_n\varphi^1-D_n\varphi^2,D^T\varphi^2-D^T\varphi_0\rangle,$$
 which has quadratic structure. Using a local frame
 $(\tau_a)_{a=1}^{n-1}$ for $T\partial D_0$, we find the components
 $\Omega_a$:
 $$\Omega_a(D\varphi^1,D\varphi^2)=[D_i(\varphi^1-\varphi^2)D_j(\varphi^1-\varphi_0)
 +D_j(\varphi^1-\varphi^2)D_i(\varphi^2-\varphi_0)]n^j\tau_a^i$$
 (summation convention, $i,j=1,\ldots n$), so $\Omega_a$ is a sum of
 terms of the form:
 $b(x)D(\varphi^1-\varphi^2)D(\varphi^3-\varphi^4)$, where
 $b(x)=n^j\tau_a^i$ and the $\varphi^I$ coincide at $t=0$.
 It is then not hard to show that:
 $$||b(x)D(\varphi^1-\varphi^2)D(\varphi^3-\varphi^4)||_{1+\alpha}\leq
 c||b||_{1+\alpha}||\varphi^1-\varphi^2||||\varphi^3-\varphi^4||T^{\alpha},$$
 with $c$ depending on the $C^1$ norms of local charts for $D_0$.
 To bound the norm $||n\otimes\tau_a||_{1+\alpha}$, note $|n||\tau_a|\leq 1$,
 $|D(n\otimes \tau_a)|\leq |Dn|+|D\tau_a|$ and $[D(n\otimes
 \tau_a)]_{x}^{(\alpha)}\leq [Dn]_x^{(\alpha)}+[D\tau_a]_x^{\alpha}$.
 Since $n=-(\beta/\beta_0)Du_0$ on $\partial D_0$
 (and $\partial D_0$ is a level set of $u_0$), we clearly have:
 $$||Dn||_{\alpha}+||D\tau_a||_{\alpha}\leq c||D^2u_0||_{\alpha}\leq
 c||u_0||.$$ We summarize the conclusion in the following lemma.

\textbf{Lemma 7.2.} Assume $\bar{\varphi},\varphi_0\in
C^{2+\alpha,1+\alpha/2}(Q^T;{\mathbb R}^n)$ have the same initial
values. Then:
 $$||\Omega(D\bar{\varphi},D\varphi_0)||_{1+\alpha}\leq
 c_0||u_0||||\bar{\varphi}-\varphi_0||^2T^{\alpha}$$
 and
 $$||\Omega(D\varphi^1,D\varphi^2)||_{1+\alpha}\leq
 c_0||u_0||(||\varphi^1-\varphi_0||+||\varphi^2-\varphi_0||)T^{\alpha}||\varphi^1-\varphi^2||,$$
 with $c_0$ depending only on the data at $t=0$. In particular, if
 $\bar{F}=[\bar{\varphi},\bar{u}]\in B_R^T$, we have:
$$||\Omega(D\bar{\varphi},D\varphi_0)||_{1+\alpha}\leq
c_0R^2T^{\alpha},$$ and for
$\bar{F}^I=[\bar{\varphi}^I,\bar{u}^I]\in B_R^T$, $I=1,2$:
$$||\Omega(D\bar{\varphi}^1,D\bar{\varphi}^2)||_{1+\alpha}\leq c_0
RT^{\alpha}||\bar{\varphi}^1-\bar{\varphi}^2||.$$
 \vspace{.3cm}

To explain the estimates for the angle condition, we write the
normal vector as a multilinear form on $DF^i$:
$$\tilde{N}(DF)=J_n(DF):=(-1)^n\sum_{i=1}^{n+1}(-1)^{i-1}(DF^1\wedge\ldots\hat{DF^i}\wedge\ldots
DF^{n+1})e_i\in{\mathbb R}^{n+1}$$ ($DF^i$ omitted in the $i^{th.}$
term of the sum), where $DF^i\in {\mathbb R}^n$ for $i=1,\ldots,
n+1$ and we identify the $n$-vector in ${\mathbb R}^n$ with a
scalar, using the standard volume form. The angle condition has the
form:
$$\beta^2|\tilde N|^2-\langle \tilde{N},e_{n+1}\rangle^2=0\mbox{ on
}\partial D_0,$$ and we set:
$$B(DF):=\beta^2|J_n(DF)|^2-\langle J_n(DF),e_{n+1}\rangle^2,$$
with linearization at $DF_0=[\mathbb{I}_n|Du_0]$:
$${\cal L}_0B[DF]=2\beta^2\langle
J_n(DF_0),DJ_n(DF_0)[DF]\rangle-2\langle
J_n(DF_0),e_{n+1}\rangle\langle DJ_n(DF_0)[DF],e_{n+1}\rangle.$$

Under the assumption $F^1=F^2$ at $t=0$, we need an estimate in
$C^{1+\alpha,\frac{1+\alpha}2}$ for:
$${\cal B}(DF^1,DF^2):=B(DF^1)-B(DF^2)-{\cal L}_0B[DF^1-DF^2]$$
$$=\beta^2(|J_n(DF^1)|^2-|J_n(DF^2)|^2-2\langle
J_n(DF_0),DJ_n(DF_0)[DF^1-DF^2]\rangle)$$
$$-(\langle J_n(DF^1),e_{n+1}^2\rangle^2-\langle
J_n(DF^2),e_{n+1}\rangle^2-2\langle J_n(DF_0),e_{n+1}\rangle \langle
DJ_n(DF_0)[DF^1-DF^2],e_{n+1}\rangle ).$$ It will suffice to
estimate the expression in the first parenthesis; the second is
analogous.\vspace{.2cm}

We need the following algebraic observation: if
$T_0=[\mathbb{I}_n|Du_0]$ and $T$ are $n\times (n+1)$ matrices, the
expression:
$$|J_n(T_0+T)|^2-|J_n(T_0)|^2-2\langle
J_n(T_0),DJ_n(T_0)[T]\rangle$$ is a linear combination (with
constant coefficients) of terms of the form:
$$u_{0i}p_{(2)}(T),\quad u_{0i}u_{0j}p_{(2)}(T),\quad p_{(2)}(T),$$
where the $p_{(2)}(T)$ are polynomials in the entries of $T$ (with
constant coefficients), with terms of degree: $2\leq deg\leq 2n$.

Thus ${\cal B}(DF^1,DF^2)$ is a linear combination (with constant
coefficients) of terms:
$$u_{0i}p_{(2)}(DF^1-DF^2),\quad u_{0i}u_{0j}p_{(2)}(DF^1-DF^2),\quad p_{(2)}(DF^1-DF^2),$$
with the $p_{(2)}$ as described; and hence is a linear combination
of terms of the form:
$$u_{0i}(F^{1j}_k-F^{2j}_k)^d,\quad
u_{0i}u_{0l}(F^{1j}_k-F^{2j}_k)^d,\quad (F^{1j}_k-F^{2j}_k)^d$$
(where $2\leq d\leq 2n,1\leq j\leq n+1,1\leq i,l,k\leq n$), which we
write symbolically as:
$${\cal B}(DF^1,DF^2)\sim \sum_{2\leq d\leq2n}b(x)(DF^1-DF^2)^d,$$
where $b(x)$ is constant or $u_{0i}(x)$ or $u_{0i}(x)u_{0j}(x).$ For
the degree $d$ terms $G^{(d)}\sim b(x)(DF^1-DF^2)^d$, it is not hard
to show the bound:
$$||G^{(d)}||_{1+\alpha}\leq
c||b||_{1+\alpha}||F^1-F^2||^dT^{\alpha},\quad 2\leq d\leq 2n.$$ We
conclude:

\textbf{Lemma 7.3.}Assume $\bar{F},F_0,F^1,F^2$ are in
 $C^{2+\alpha,1+\alpha/2}(Q^T;{\mathbb R}^{n+1})$ and have the same
 initial values. Then:
$$||{\cal B}(D\bar{F},DF_0)||_{1+\alpha}\leq
c(1+||u_0||^2)(1+||\bar{F}-F_0||^{2n-2})T^{\alpha}||\bar{F}-F_0||^2.$$
$$||{\cal B}(DF^1,DF^2)||_{1+\alpha}\leq
c(1+||u_0||^2)(1+||F^1-F^2||^{2n-2})T^{\alpha}||F^1-F^2||^2,$$ with
$c$ depending only on $F_0$. In particular, if $\bar{F}\in B_R^T$:
$$||{\cal B}(D\bar{F},DF_0)||_{1+\alpha}\leq c_0 R^2T^{\alpha},$$
and if $\bar{F}^1,\bar{F}^2\in B_R^T$:
$$||{\cal B}(D\bar{F}^1,D\bar{F}^2)||_{1+\alpha}\leq
c_0T^{\alpha}||\bar{F}^1-\bar{F}^2||,$$ with $c_0$ depending only on
$F_0$.
 \vspace{.3cm}

\textbf{8.} \textbf{Local existence.}

Let $D_0\subset {\mathbb R^n}$ be a uniformly $C^{2+\alpha}$ domain,
not necessarily bounded or connected (note: we define our norms as
the sum of the norms on each connected component).\vspace{.2cm}

Given a $C^{3+\alpha}$ graph $\Sigma_0$ over $D_0$ satisfying the
contact and angle conditions, let $\varphi_0\in Diff^{2+\alpha}$ be
a diffeomorphism given by lemma 4.1 (with the 1-jet of the identity
at $\partial D_0$ and 2-jet determined by the mean curvature of
$\Sigma_0$ at $\partial D_0$). Then find $u_0\in C^{2+\alpha}(D_0)$
so that $F_0=[\varphi_0,u_0]$ parametrizes $\Sigma_0$ over $D_0$.

(Precisely, if $[z,\tilde{u_0}(z)]$ parametrizes $\Sigma_0$ as a
graph, and $\varphi_0$ is given by lemma 4.1, let
$u_0=\tilde{u_0}\circ \varphi_0$; so $u_0\in C^{2+\alpha}$.)
 \vspace{.2cm}

We obtained in section 7 all the estimates needed for a fixed-point
argument in the set:
$$B_R^T=\{F\in C^{2+\alpha,1+\alpha/2}(Q^T,{\mathbb
R}^{n+1});||F-F_0||<R, F|_{t=0}=F_0\}.$$

Choose $R<1$ and $T_0<1$ small enough (depending only on $F_0$) so
that, for $F\in B_R^{T_0}$, $F(t)=[\varphi(t),u(t)]$ defines an
embedding of $D_0$, with $\varphi(t)$ a diffeomorphism onto its
image $D(t)$. Let $K\subset Imm(\mathbb{R}^n,\mathbb
 {R}^{n+1})$ be a compact set containing $DF(z)$ for all $F\in B_R,
 z\in Q^{T_0}$. Now consider $T<T_0$.
\vspace{.2cm}

Given $\bar{F}\in B_R^T$, solve the linear system (with initial data
$F_0$) to obtain $F\in C^{2+\alpha,1+\alpha/2}(Q^T)$. (This is
possible since the complementarity and compatibility conditions hold
for the linear system.) This defines a map ${\mathbb F}:{\bar
F}\mapsto F$.

From linear parabolic theory (e.g. \cite{EidelmanZhitarasu}, thm
VI.21]):
$$||F-F_0||\leq M(||{\cal F}(\bar{F},F_0)||_{\alpha}+||{\cal
B}(D\bar{F},DF_0)||_{1+\alpha}+||\Omega(D\bar
{\varphi},D\varphi_0)||_{1+\alpha}),$$ where $M>0$ depends on the
$C^{\alpha,\alpha/2}$ norm of the coefficients of the linear system,
that is, ultimately on $||F_0||$. \vspace{.2cm}

From lemmas 7.2-7.4 in section 7, it follows that:
$$||F-F_0||\leq Mc_0(RT^{\alpha/2}+R^2T^{\alpha})<R$$
provided $T$ is chosen small enough (depending only on $F_0$.) Thus
$\mathbb{F}$ maps $B_R^T$ to itself.\vspace{.2cm}

Similarly, if ${\mathbb F}(\bar{F}^i)=F^i$ for $i=1,2$, standard
estimates for the linear system solved by $F^1-F^2$ give:
$$||F^1-F^2||\leq M(||{\cal F}(\bar{F}^1,\bar{F}^2)||_{\alpha}+||{\cal
B}(D\bar{F}^1,D\bar{F}^2)||_{1+\alpha}+||\Omega(D
{\bar{\varphi}}^1,D\bar{\varphi}^2)||_{1+\alpha})$$

Again the estimates in lemmas 7.2-7.4 imply:
$$||F^1-F^2||\leq Mc_0(T^{\alpha/2}+T^{\alpha})||\bar{F}^1-\bar{F}^2||<\frac 12||\bar{F}^1-\bar{F}^2||,$$
assuming $T$ is small enough (depending only on $F_0$). This
concludes the argument for local existence. \vspace{.2cm}

\textbf{Theorem 8.1.} Let $\Sigma_0\subset {\mathbb R}^{n+1}$ be a
$C^{3+\alpha}$ graph over $D_0\subset {\mathbb R}^n$ satisfying the
contact and angle conditions at $\partial D_0$ ($\Sigma_0$ may be
unbounded or not connected). There exists a parametrization
$F_0=[\varphi_0,u_0]\in C^{2+\alpha}(D_0)$ of $\Sigma_0$, $T>0$
depending only on $F_0$ and a unique solution $F\in
C^{2+\alpha,1+\alpha/2}(Q^T;{\mathbb R}^{n+1})$ of the system:
$$\left \{ \begin{array}{l}
\partial_t F-g^{ij}(DF)\partial_i\partial_jF=0,\\
u_{|\partial D_0}=0, \quad N^{n+1}(D\varphi,Du)_{|\partial
D_0}=\beta,\end{array}\right .$$ with initial data $F_0$. For each
$t\in [0,T)$, $F(t)$ is a $C^{2+\alpha}$ embedding parametrizing a
surface $\Sigma_t$ which satisfies the contact and angle conditions
and moves by mean curvature. In addition, $F(t)$ satisfies the
orthogonality condition at $\partial D_0$.

The hypersurfaces $\Sigma_t$ are graphs. For each $t\in [0,T)$,
$\varphi(t):D_0\rightarrow D(t)$ is a diffeomorphism and
$\Sigma_t=graph(w(t))$, for $w(t):D(t)\rightarrow {\mathbb R}$ given
by $w(t)=u(t)\circ \varphi^{-1}(t)$. (We have $w(t)\in
C^{2+\alpha^2}(D(t))$, `less regular' than $u(t)$ or $\varphi(t)$.)
$D(t)$ is a uniformly $C^{2+\alpha}$ domain.\vspace{.3cm}

\emph{Remark 8.1.} This theorem does not address geometric
uniqueness of the motion, given $\Sigma_0$. It only asserts
uniqueness for solutions of the parametrized flow (including the
orthogonality boundary condition) in the given regularity
class.\vspace{.3cm}

\textbf{9. Rotational symmetry.} In this section we record the
equations for two rotationally symmetric instances of the problem:
(i) $D_0$ and $D(t)$ are disks, and $u>0$ (`lens' case); (ii) $D_0$
and $D(t)$ are complements of disks in ${\mathbb R}^n$ (`exterior'
case). For simplicity we restrict to $n=2$.

Let $F(r)=[\varphi(r),u(r)]$ parametrize a hypersurface $\Sigma$,
where $\varphi(r)=\phi(r)e_r$ is a diffeomorphism onto its image.
Here $e_r,e_{\theta}$ are orthonormal vectors, outward normal (resp.
counterclockwise tangent) to the circles $r$=const. The unit upward
normal vector and mean curvature are:
$$N=\frac{[-u_re_r,\phi_r]}{\sqrt{u_r^2+\phi_r^2}},$$
$$H=\frac 1{(\phi_r^2+u_r^2)^{3/2}}(\phi_r{\cal
M}(\phi_r,u_r)[D^2u]-\langle u_re_r,\vec{\cal
M}(\phi_r,u_r[D^2\varphi]\rangle),$$ where:
$${\cal
M}(\phi_r,u_r)[D^2u]=u_{rr}+(\phi_r^2+u_r^2)\frac{u_r\phi_r}{\phi^2},$$
$$\vec{\cal
M}(\phi_r,u_r)[D^2\varphi]=[\phi_{rr}+(\phi_r^2+u_r^2)(\frac{r\phi_r}{\phi^2}-\frac
1{\phi})]e_r.$$ Simplifying:
$$H=\frac
1{(\phi_r^2+u_r^2)^{3/2}}[\phi_ru_{rr}-u_r\phi_{rr}+(\phi_r^2+u_r^2)\frac{u_r}{\phi}].$$

Now consider the time-dependent case $F(r,t)=[\phi(r,t)e_r,u(r,t)]$.
From the above expressions, one finds easily that the equation
$\langle\partial_tF,N\rangle=H$ takes the form:
$$\phi_r(u_t-\frac 1{\phi_r^2+u_r^2}{\cal
M}(\phi_r,u_r)[D^2u])=u_r\langle e_r,\varphi_t-\frac
1{\phi_r^2+u_r^2}\vec{\cal M}(\phi_r,u_r)[D^2\varphi]\rangle.$$

In `split gauge', we consider the system:
$$u_t-\frac 1{\phi_r^2+u_r^2}{\cal
M}(\phi_r,u_r)[D^2u]=0,$$
$$\varphi_t-\frac
1{\phi_r^2+u_r^2}\vec{\cal M}(\phi_r,u_r)[D^2\varphi]=0.$$

Note that $\phi(r,t)=r$ solves the $\phi$ equation, and that in this
case the $u$ equation becomes:
$$w_t-\frac{w_{rr}}{1+w_r^2}-\frac{w_r}r=0.$$ This can be compared
with the equation for curve networks:
$$w_t-\frac{w_{xx}}{1+w_x^2}=0.$$\vspace{.2cm}

 The boundary conditions are easily stated (we assume $D_0$ is the unit disk or its
 complement).

 The `contact condition' at $r=1$
is $u=0$. For the `angle condition' at $r=1$, we find:
$$u_r^2=\frac{\beta_0^2}{\beta^2}\phi_r^2,\quad
\beta_0:=\sqrt{1-\beta^2}.$$ Assuming $\phi_r>0$ at $r=1$, this
resolves as:
$$\beta u_r+\beta_0\phi_r=0\mbox{ at }r=1\mbox{ (lens case)};$$
$$\beta u_r-\beta_0\phi_r=0\mbox{ at }r=1\mbox{ (exterior case)}.$$
(For lenses, one also has at $r=0$: $u_r=0$ and $\phi_r=1$.) Thus in
both cases one can work with \emph{linear} Dirichlet/Neumann-type
boundary conditions.\vspace{.2cm}

One reason to consider the exterior case is that (unlike the lens
case) it admits stationary solutions. Geometrically, one just has to
consider one-half of a catenoid, truncated at an appropriate height.
For example, for 120 degree junctions the equation for stationary
solutions:
$$\frac{u_{rr}}{1+u_r^2}+\frac{u_r}r=0\mbox{ in }\{r>1\},$$
$${u_r}_{|r=1}=\sqrt{3},\quad u_{|r=1}=0$$
admits the explicit solution:
$$u(r)=\frac{\sqrt{3}}2(\ln(2r+\sqrt{4r^2-3})-\ln 3),\quad
r>\sqrt{3}/2.$$\vspace{.3cm}

\emph{Problem.} It would be interesting to consider the nonlinear
dynamical stability of this solution (even linear stability is yet
to be considered.) One may even work with bounded domains, by
introducing a fixed boundary at some $R>1$, intersecting the surface
orthogonally (see Section 10.)\vspace{.3cm}

\textbf{10. Fixed supporting hypersurfaces.} Extending the local
existence theorem to the case of hypersurfaces intersecting a fixed
hypersurface $\cal S$ orthogonally presents no essential difficulty.
The case of vertical support surface leads directly to graph
evolution with a standard Neumann condition on a fixed boundary; we
consider the complementary case where $\cal S$ is a graph. Let
${\cal S}\subset {\mathbb R}^{n+1}$ be a $C^{4}$ embedded
hypersurface (not necessarily connected), the graph over ${\cal
D}\subset {\mathbb R}^n$ of $B\in C^4({\cal D})$, oriented by the
upward unit normal:
$$\nu(y):=\frac 1{v_B}\tilde{\nu}(y),\quad
\tilde{\nu}(y):=[-DB(y),1]\in {\mathbb R}^n\times {\mathbb R}, \quad
v_B:=\sqrt{1+|DB(y)|^2}.$$ $\nu$ is assumed to be nowhere vertical
in ${\cal D}$ ($DB\neq 0$).  To state the problem in the graph
parametrization, we consider a time-dependent domain $D(t)\subset
{\mathbb R}^n$ with boundary consisting of two components
$\partial_1D(t)$ and $\partial_2D(t)$, both moving. The hypersurface
$\Sigma_t$ is the graph of $w(\cdot, t)$ over $D(t)$, solving the
parabolic equation:
$$w_t-g^{ij}(Dw)w_{ij}=0\mbox{ in }E:=\bigcup_{t\in [0,T]}D(t)\times
\{t\}\in {\mathbb R}^{n+1}\times [0,T],$$ with boundary conditions:
$$w(\cdot,t)_{|\partial_1D(t)}=0,\quad
\sqrt{1+|Dw|^2}_{|\partial_1D(t)}=1/\beta$$ (as before), and on
$\partial_2D(t)$:
$$w=B, \quad \nabla w\cdot \nabla B=-1.$$
 (The first-order condition on $\partial_2D(t)$ is equivalent to
$\langle \nu,N\rangle=0$).\vspace{.2cm}

Differentiating in $t$ the boundary condition $w=B$ leads easily to
an equation for the normal velocity of the interface
$\Gamma(t)=\partial_2D(t)$:
$$\dot{\Gamma_n}=\frac{vH}{B_n-w_n}.$$\vspace{.2cm}
Note that $w_n$ at $\partial_2D(t)$ can be computed from $B_n$,
since:
$$-1=\nabla w\cdot \nabla B=w_nB_n+|\nabla^TB|^2;$$
in particular neither $B_n$ nor $w_n$ can vanish (so both have
constant sign on connected components of $\partial_2D$), and one
easily computes: $w_n-B_n=-v_B^2/B_n$.\vspace{.2cm}

Let $\Lambda=\Sigma\cap {\cal S}$ be the $(n-1)$-manifold of
intersection, the graph of $w$ (or $B$) over $\partial_2D$. Given
the graph parametrizations of $\Sigma$ and $\cal S$:
$$G(y)=[y,w(y)],\quad {\mathbb B}(y)=[y,B(y)],\quad y\in
\partial_2D,$$
and $\tau \in T\partial_2 D$, we have the tangent vectors:
$$G_n:=[n,w_n]\in T\Sigma,\quad G_B:=[\nabla B,-1]=-v_B\nu\in
T\Sigma,\quad G_{\tau}:=[\tau, \nabla w\cdot \tau]\in T\Lambda,$$
and the second fundamental forms of $\Sigma$ and ${\cal S}$(for
$e\in {\mathbb R}^n$ arbitrary):
$$A(dGe,dGe)=\frac 1{v}d^2w(e,e),\quad {\cal A}(d{\mathbb
B}e,d{\mathbb B}e)=\frac 1{v_B}d^2B(e,e).$$ From $\langle
\nu,N\rangle=0$ at $\partial_2D$, it follows easily that (cp.
\cite{Stahl}):
$$A(G_\tau,\nu)=-{\cal A}(G_\tau,N),
\quad \tau\in T\partial D.$$

For the remainder of this section, we concentrate on the boundary
conditions at $\partial_2D_0$, and denote this boundary component
simply by $\partial D_0$. To establish short-time existence, we
consider as before the parametrized flow:
$$F_t-tr_gd^2F=0,\quad g=g(dF),\quad F=[\varphi,u].$$
The contact and angle boundary conditions are:
$$u_{|\partial D_0}=B\circ \varphi_{|\partial D_0},\quad \langle
N,\nu\circ \varphi\rangle_{|\partial D_0}=0.$$ Again we have two
scalar boundary conditions for $n+1$ components. Here the solution
is easier than at the junction. With the notation
$F_n=dFn=[\varphi_n,u_n]$, we replace the angle condition by the
`vector Neumann condition':
$$F_n\perp T{\cal S},\mbox{ or }F_n=-\alpha v_B\nu\mbox{ on
}\partial D_0,$$ where $\alpha:\partial D_0\rightarrow {\mathbb R}$,
or equivalently (since this leads to $\alpha=-u_n$):
$$\varphi_n=-u_n(\nabla B\circ \varphi) \mbox{ on }\partial_2D_0.$$
Clearly the Neumann condition implies the angle condition $\langle
N,\nu\circ \varphi\rangle=0$, but not conversely. This linear
Neumann-type condition can easily be incorporated into the
fixed-point existence scheme described earlier. \vspace{.2cm}

There is one issue to consider: the 0 and 1st-order compatibility
conditions must hold at $\partial D_0$, at $t=0$. The initial
hypersurface $\Sigma_0$ uniquely determines $w_0$ and $D_0\subset
\mathbb{R}^n$ (satisfying $w_0=B$ and $\nabla w_0\cdot \nabla B=-1$
on $\partial D_0$), and then once $\varphi_0\in Diff(D_0)$ is fixed,
$u_0=w_0\circ \varphi_0$ is also determined. We may assume:
$${\varphi_0}=id, \quad \varphi_{0n}=\nabla B\mbox{ on }\partial
D_0,$$ so:$$u_{0n}=\nabla w_0\cdot \varphi_{0n}=\nabla w_0\cdot
\nabla B=-1\mbox{ on }\partial D_0,$$ and then the Neumann condition
${F_{0n}}_{|\partial_2D_0}=-v_B\nu$ holds at $t=0$, on $\partial
D_0$.\vspace{.2cm}

The first-order compatibility condition is:
$$tr_gd^2u_0=u_t=\nabla B\cdot \varphi_t=\nabla B\cdot
tr_gd^2\varphi_0\mbox{ on }\partial D_0,$$ or equivalently:
$$tr_g\langle \nu,d^2F_0\rangle=0\mbox{ on }\partial D_0.$$
(This is not a mean curvature condition; the mean curvature of
$\Sigma_0$ is $H=tr_g\langle N,d^2F_0\rangle.$) \vspace{.2cm}

From now on we omit the subscript $0$, but continue to discuss
compatibility at $t=0$. First observe that the Neumann condition
leads to a splitting of the induced metric. Given $\tau \in
T\partial D_0$, let $F_{\tau}=dF\tau\in T\Lambda$. Then (recalling
$u_n=-1$ on $\partial D_0$):
$$\langle F_{\tau},F_n\rangle=\langle[\tau, dB
\tau],[\varphi_n,u_n]\rangle=\nabla B\cdot \tau-\nabla B\cdot
\tau=0.$$ Thus we have:
$$tr_g\langle \nu,d^2F\rangle=g^{ab}\langle
\nu,d^2F(\tau_a,\tau_b)\rangle +g^{nn}\langle
\nu,d^2F(F_n,F_n)\rangle,$$ for a local basis
$\{T_a=dF\tau_a\}_{a=1}^{n-1}$ of $T\Lambda$, with $g_{ab}=\langle
T_a,T_b\rangle$ and $g_{nn}=|F_n|^2=v_B^2$.\vspace{.2cm}

Differentiating in $n$ the condition $u_n=\nabla w\cdot \varphi_n$ (
assuming, as usual, $n$ extended to a tubular neighborhood $\cal N$
of ${\partial D}_0$ as a self-parallel vector field), we find:
$$u_{nn}=d^2w(n,\nabla B)+ \nabla w\cdot d^2\varphi(n,n).$$
This is used to compute:
$$\langle \nu,d^2F(n,n)\rangle =\frac 1{v_B}[u_{nn}-\nabla B\cdot
d^2\varphi(n,n)]$$
$$=\frac 1{v_B}[d^2w(n,\nabla B)+(\nabla w-\nabla
B)\cdot d^2\varphi(n,n)]$$
$$=-vA(G_n,\nu)+\frac 1{v_B}(w_n-B_n)
n\cdot d^2\varphi(n,n).$$ Bearing in mind the expression for
$w_n-B_n$ found earlier, the compatibility condition may be stated
in the form:
$$\frac{v_B}{B_n}
n \cdot d^2\varphi(n,n)=-vA(G_n,\nu)+g^{ab}\langle
d^2F(\tau_a,\tau_b),\nu\rangle.$$\vspace{.2cm}

We are now in the same situation as in section 4: given the 1-jet of
$\varphi_0$ on $\partial D_0$, we extend $\varphi_0$ to a tubular
neighborhood $\cal N$ of $\partial D_0$ (and then to all of $D_0$),
so that $ n\cdot d^2\varphi(n,n)$ has on $\partial D_0$ the value
dictated by the compatibility condition (using Lemma 4.1(ii)). We
just need to verify that the right-hand side of the above expression
depends only on $\Sigma_0$, $\cal S$, and the 1-jet of $\varphi_0$
over $\partial D_0$. Clearly only the term $g^{ab}\langle
\nu,d^2F(\tau_a,\tau_b)\rangle$ is potentially an issue.
\vspace{.2cm}

Fix $p\in \partial D_0$, and let $\{\tau_a\}$ be an orthonormal
frame for $T\partial D_0$ near $p$, parallel at $p$ for the
connection induced on $\partial D_0$ from ${\mathbb R}^n$. If $\cal
K$ denotes the second fundamental form of $\partial D_0$ in
${\mathbb R}^n$, we have:
$$\tau_a(\tau_b)={\cal K}(\tau_a,\tau_b)n\quad (\mbox{ at }p)$$
(on the left-hand-side, $\tau_b$ is regarded as a vector-valued
function in ${\mathbb R}^n$). Still computing at $p$, this implies:
$$d^2F(\tau_a,\tau_b)=\tau_a(dF\tau_b)-dF(\tau_a(\tau_b))$$
$$=\tau_a(d{\mathbb B}\tau_b)-{\cal K}(\tau_a,\tau_b)F_n$$
$$=d^2{\mathbb B}(\tau_a,\tau_b)+{\cal K}(\tau_a,\tau_b){\mathbb
B}_n-{\cal K}(\tau_a,\tau_n)F_n,$$ where $F_n=-v\nu$ and ${\mathbb
B}_n=d{\mathbb B}n\in T{\cal S}$. Hence:
$$\langle \nu,d^2F(\tau_a,\tau_b)\rangle=\langle \nu,d^2{\mathbb
B}(\tau_a,\tau_b)\rangle +v{\cal K}(\tau_a,\tau_b)={\cal
A}(T_a,T_b)+v{\cal K}(\tau_a,\tau_b).$$ This clearly depends only on
$\cal S$ and on $\Sigma_0$. \vspace{.2cm} We summarize the
discussion in a lemma.\vspace{.2cm}

\textbf{Lemma 10.1} Let $\Sigma_0=graph(w_0)$ be a $C^3$ graph over
$D_0\subset {\mathbb R}^n$ (a uniformly $C^3$ domain), intersecting
a fixed hypersurface ${\cal S}=graph(B)$ over $\partial D_0$.
Consider the parametrized mean curvature motion with Neumann
boundary condition:
$$F\in C^{2,1}(D_0\times [0,T])\rightarrow {\mathbb R}^{n+1},\quad
F=[\varphi,u]$$
$$F_t-tr_gd^2F=0,\quad g=g(dF),\quad u\circ \varphi=B\mbox{ and }F_n\perp T{\cal S}\mbox{
on }\partial D_0.$$ Then $\varphi_0\in Diff(D_0)$ can be chosen so
that (with $u_0=w_0\circ \varphi_0$) the initial data
$F_0=[\varphi_0,u_0]$ satisfies the order zero and the first-order
compatibility conditions at $t=0$ and $\partial D_0$:
$$\varphi_{0n}=-u_{0n}(\nabla B\circ {\varphi_0}),\quad \langle
\nu\circ \varphi_0,tr_{g_0}d^2F_0\rangle=0.$$

\emph{Remark 10.1.} Differentiating $dw\tau_a=dB\tau_a$ along
$\tau_b$, we find:
$$d^2w(\tau_a,\tau_b)-d^2B(\tau_a,\tau_b)=(w_n-B_n){\cal
K}(\tau_a,\tau_b)$$ (reminding us that, although $w\equiv B$ on
$\partial D_0$, the tangential components of their Hessians do not
coincide.) From this follows the expression for ${\cal K}$ in terms
of $A$ and $\cal A$:
$${\cal K}(\tau_a,\tau_b)=\frac
1{w_n-B_n}[vA(T_a,T_b)-v_B{\cal A}(T_a,T_b)].$$ It is also easy to
express the corresponding traces in terms of the mean curvatures
$H^{\Lambda}$ and ${\cal H}^{\Lambda}$ of $\Lambda$ in $\Sigma$ and
${\cal S}$:
$$H^{\Lambda}=\frac v{v_B}g^{ab}A(T_a,T_b),\quad {\cal
H}^{\Lambda}=\frac{v_B}vg^{ab}{\cal A}(T_a,T_b).$$\vspace{.3cm}

\textbf{11. A continuation criterion.} Once local existence has been
established, it is easier to obtain geometric estimates (in
particular using the maximum principle) for the solution in the
graph parametrization. (From this point on , we focus on the `lens'
case, without fixed support hypersurfaces.)

For a time interval $I=(t_0,t_1)\subset [0,T]$ set:
$$E=\{z=(y,t)\in {\mathbb R}^n\times I;y\in D(t)\},\quad S=\{(y,t);t\in I,y\in \partial D(t)\}.$$

Let $w$ be a solution in $E$ of:
$$w_t-g^{ij}(Dw)D^2_{i,j}w=0$$
with boundary conditions on $S$:
$$w=0,\quad D_nw=\beta_0v,\quad \beta_0:=\sqrt{1-\beta^2}.$$
For the remainder of the paper we assume $D(t)$ is \emph{bounded},
for each $t\in I$. $n$ denotes the ($t$-dependent) \emph{inner} unit
normal at $\partial D(t)$, extended to a $C^{2,1}$ unit vector field
in a tubular neighborhood of $D(t)$ so that $D_nn=0$.\vspace{.2cm}

Denote by $L$ the operator
$L=\partial_t-g^{ij}(Dw)\partial_i\partial_j$, so $Lw=0$ in $E$. The
following height bound is immediate.

\textbf{Lemma 11.1.} Assume $0<w_0<M$ in $D(t_0)$. (If there is a
support surface $\cal S$, we assume $B(y)>0$ in ${\cal D}$ and
$M<\sup_{\cal D}B$.) Then $0<w<M$ in $\bar{E}$.

\emph{Proof.} Follows from the maximum principle applied to $L$,
since $0\leq w\leq M$ holds on the parabolic boundary $\partial_pE$.
\vspace{.2cm}

It is well-known that the function $v=\sqrt{1+|Dw|^2}$ solves the
evolution equation (assuming $Dw\in C^{2,1}(\bar{E})$, see e.g.
\cite{Guan}):

$$Lv+\frac 2vg^{ij}v_iw_j=-v|A|^2_g.$$
From the maximum principle, we have the following global bound on
$v$ (equivalently, on $|Dw|$).\vspace{.2cm}

\textbf{Lemma 11.2} Assume $w$ is a solution with $Dw\in
C^{2,1}(\bar{E})$. Then we have on $\bar{E}$:
$$v(z)\leq \max\{sup_{D(t_0)}v(x,t_0),\frac 1{\beta}\}.$$

\emph{Proof.} By the maximum principle,
$\max_{\bar{E}}v=\max_{\partial_pE}v$. Note $v_{|S}\equiv \frac
1{\beta}.$\vspace{.2cm}

It follows from this lemma that $g_{ij}(t)$ is uniformly equivalent
to the euclidean metric in $D(t)$: if $v\leq \bar{v}$ in $\bar{E}$,
and $X$ is a vector field in $D(t)$:
$$|X|_e^2\leq |X|_g^2=g_{ij}X^iX^j=|X|^2_e+(X\cdot Dw)^2\leq
|X|^2_e(1+|Dw|^2)\leq \bar{v}^2|X|_e^2.$$ Also, if
$\omega:=v^{-1}Dw$:
$$|\omega|^2_e=\frac{|Dw|^2_e}{v^2}=1-\frac 1{v^2}\leq
1-\frac1{\bar{v}^2}.$$ This equivalence of norms clearly extends to
tensors, in particular to $h$:
$$\frac 1{c_n}|h|^2_e\leq |h|_g^2\leq c_n|h|^2_e,$$
where (throughout this section) $c_n$ denotes a constant depending
only on $n$ and $\bar{v}$. More generally, defining:
$$|\partial h|^2_e:=\sum_{i,j,k}(\partial_kh_{ij})^2,\quad |\partial
h|^2_g:=\sum_{i,j,k,l}g^{kl}(\partial_k h_{ij})(\partial_lh_{ij}),$$
we have, for each pair $i,j$:
$$\sum_{k}(\partial_kh_{ij})^2=\sum_{k,l}\delta^{kl}\partial_kh_{ij}\partial_lh_{ij}
=\sum_{k,l}(g^{kl}+\omega^k\omega^l)\partial_kh_{ij}\partial_lh_{ij}$$
$$=(d_{\omega}h_{ij})^2+\sum_{k,l}g^{kl}\partial_kh_{ij}\partial_lh_{ij}\leq
\sum_{k,l}g^{kl}\partial_kh_{ij}\partial_lh_{ij}+(1-\frac
1{\bar{v}}^2)\sum_{k}(\partial_kh_{ij})^2,$$ and hence, adding over
$i,j$ we have:
$$|\partial h|^2_g\leq |\partial h|^2_e\leq \bar{v}^2|\partial
h|_g^2.$$ The same argument works for second derivatives. The norms
defined by:
$$|\partial^2h|_e^2:=\sum_{i,j,k,m}(\partial_m\partial_kh_{ij})^2,\quad
|\partial^2h|^2_g:=\sum_{i,j,k,m,n}g^{mn}(\partial_m\partial_kh_{ij})(\partial_n\partial_kh_{ij})$$
are uniformly equivalent in $E$:
$$|\partial^2h|^2_g\leq |\partial^2h|^2_e\leq
\bar{v}^2|\partial^2h|^2_g.$$ \vspace{.2cm}

The point is that these euclidean norms satisfy easily computed
evolution equations. Using the results in Appendix 2, it is
straightforward to see that:
$$L[|h|^2_e]=-2|\partial h|_g^2+2\sum_{i,j}C_{ij}h_{ij},\quad
C_{ij}:=L[h_{ij}],$$
$$L[|\partial h|^2_e]=-2|\partial^2
h|^2_g+2\sum_{i,j,k}(\partial_kC_{ij})(\partial_kh_{ij})
+2\sum_{i,j,k,m,n}(\partial_kg^{mn})(\partial_m\partial_n
h_{ij})(\partial_kh_{ij}).$$ In symbolic notation, we have:
$$C_{ij}\sim h*h*h,\quad \partial_kC_{ij}\sim (\partial h)*h*h,$$
which combined with the previous remarks implies:
$$\sum_{ij}(C_{ij})^2\leq c_n|h|_g^6,\quad
\sum_{i,j,k}(\partial_kC_{ij})^2\leq c_n|\partial h|^2_g|h|_g^4,$$
for a constant $c_n$ as above. In addition, from (see Appendix 2):
$$\partial_kg^{mn}=h_k^m\omega^n+h_k^n\omega^m\mbox{ and
}|g^{ij}|\leq 2,|\omega^n|<1,$$ we have $|\partial_kg^{mn}|\leq
c_n|h|_g$ for each $m,n,k$. We conclude:
$$L[|h|^2_e]\leq -2|\partial h|_g^2+c_n |h|_g^4,$$
$$L[|\partial h|^2_e]\leq -2|\partial^2h|_g^2+c_n |\partial
h|_g^2|h|_g^2+c_n|h|_g|\partial^2h|_g|\partial h|_g.$$ These
differential inequalities imply the continuation criterion given in
Proposition 11.3.\vspace{.2cm}

Recall that for mean curvature flow (or mean curvature motion) of
graphs, interior estimates for $v$ imply interior estimates for $A$
and its covariant derivatives $\nabla^mA$ of any order (see
\cite{EckerHuisken} or \cite{EckerBook}). In the following
`continuation criterion', global bounds are needed. \vspace{.2cm}

\textbf{Proposition 11.3.} Assuming $T_{max}$ is finite, let
$w:E^{T_{max}}\rightarrow {\mathbb R}$ be a maximal solution,
defined for $t\in [0,T_{max})$. Then:
$$\limsup_{t\rightarrow T}(\sup_{(y,t)\in E}|h|_g+\sup_{y\in \partial D(t)}|\nabla
h|_g(y,t))=\infty.$$

\emph{Proof.} By contradiction, assume we have bounds in $[t_0,t_1]$
(for $t_1$ arbitrarily close to $T_{max}$):
$$\sup_{z\in E}|h|_g\leq a_0,\quad \sup_{z\in S\cup D(0)}|\nabla
h|_g\leq b_0.$$ For $\alpha>0$ to be chosen (small), define the
function on $E$:
$$f(x,t)=\alpha|\partial h|_e^2+|h|_e^2.$$
Then, for any $\eta>0$:
$$L[f]\leq -2\alpha |\partial^2h|_g^2-2|\partial
h|_g^2+c_n(a_0^4+\alpha a_0^2|\partial h|_g^2+\alpha
a_0|\partial^2h|_g|\partial h|_g)$$
$$\leq -2\alpha |\partial^2 h|^2_g-2|\partial h|_g^2+c_na_0\alpha
\eta|\partial^2h|_g^2+(c_n\alpha
a_0^2+\frac{c_na_0\alpha}{\eta})|\partial h|_g^2+c_na_0^4.$$
Choosing $\eta$ so that $c_na_0\eta\leq 1$, then $\alpha$ so that
$c_n\alpha a_0^2+\frac{c_na_0\alpha}{\eta}<1$, we ensure that:
$$L[f-c_na_0^4t]\leq 0$$
in $E$. By the maximum principle:
$$\alpha \sup_E |\partial h|^2_e\leq \sup_Ef\leq
\sup_{\partial_pE}f+c_na_0^4T\leq c_n(a_0^2+\alpha b_0^2+a_0^4T).$$

This implies a uniform $C^3(\bar{D}(t))$ bound for $w$ in $\bar{E}$,
and hence (by linear parabolic theory, given the uniform bound on
$|Dw|$ from lemma 10.2) a $C^{3+\alpha}$ bound for some $0<\alpha
<1$. So we can apply the local existence theorem with initial data
$\Sigma_{t_1}$, to continue the solution for a time depending only
on bounds at $t_0$, contradicting the maximality of
$T_{max}$.\vspace{.2cm}

Lemma 12.2 (in the next section) implies the conclusion can be
strengthened: only a uniform bound on tangential covariant
derivatives of the second fundamental form $\cal K$ of the moving
boundary (in $\mathbb{R}^n$) is needed:\vspace{.2cm}

\textbf{Proposition 11.4.} Assuming $T_{max}$ is finite:
$$\limsup_{t\rightarrow T_{max}}[\sup_{(y,t)\in E}|h|_g+\sup_{y\in \partial
D(t)}|\nabla_{\tau} {\cal K}|_g(y,t)]=\infty.$$\vspace{.3cm}

 It is possible to strengthen this further and show that:
$$\limsup_{t\rightarrow T_{max}}[\sup_{y\in \partial D(t)}|h|_g(y)+\sup_{y\in \partial
D(t)}|\nabla_{\tau} {\cal K}|_g(y,t)]=\infty.$$ That is, the
supremum of $|h|_g$ on the moving boundary controls its value in the
interior. The reason is that we already have a bound on $\sup_E v$;
as remarked earlier, it is a well known-fact for mean curvature flow
of graphs that this implies interior bounds for the second
fundamental form and its covariant derivatives (\cite{EckerHuisken},
\cite{EckerBook}). In the next lemma we describe a global argument
for mean curvature motion of graphs with moving
boundaries.\vspace{.2cm}

\textbf{Proposition 11.5.} Let $w:E\rightarrow \mathbb{R}$ be a
solution of graph m.c.m in a spacetime domain $E\subset
\mathbb{R}^n\times [0,T]$, where $T<\infty$. Assume the first
derivative bound $v(x,t)\leq \bar{v}$ holds globally in $\bar{E}$.
Then if the bound $|h|_g\leq h_0$ holds on the parabolic boundary
$\partial_pE$, we also have the global bound:
$$|h|_g\leq a_0\quad \mbox{ in }\bar{E},$$
for a constant $a_0$ depending only on $n,\bar{v},h_0,T$ and the
initial data of $w$. \vspace{.2cm}

\emph{Proof.} The idea is to consider a function in $E$ of the form:
$$\varphi=Av^p+|h|^2_gv^p+B|h|^2_g,$$
where $A$, $B$ and $p$ are positive constants. We claim it is
possible to choose these constants, and also $C>0$ (all depending
only on $n$ and $\bar{v}$) so that: $$L[\varphi]\leq C\mbox{ in
}E.$$ Thus $L[\varphi-Ct]\leq 0$ in $E$, and hence by the maximum
principle:
$$\sup_EB|h|^2_g\leq \sup_E\varphi\leq
sup_{\partial_pE}\varphi+CT,$$ which clearly implies the bound
claimed in the proposition. \vspace{.2cm}

The proof that $\varphi$ as above exists is (of course) based on the
evolution equations for $|h|_g$ and $v$ (see Appendix 2), which
imply (for constants $c_n$, $d_n$ depending only on $n$):
$$L[|h|^2_g]\leq -2|\nabla h|^2_g+c_n|h|^4_g,$$
$$L[v^p]=-pv^p|h|^2_g-p(p-1)v^{p-2}|\partial
v|_g^2-2pv^{p-2}g^{ij}v_iw_j.$$ Here $|\partial
v|_g^2:=g^{kl}v_kv_l$, and we have the bounds:
$$L[v^p]\leq -pv^p|h|^2_g-p(p-1)v^{p-2}|\partial
v|_g^2+d_npv^{p-1}|\partial v|_g,$$
$$L[v^p]\leq -pv^p|h|_g^2-[p(p-1)-\frac 14]|\partial
v|_g^2+d_n^2p^2(\bar{v})^{2(p-1)}.$$

The main term in $L[\varphi]$ is:
$$L[|h|^2_gv^p]=L[|h|^2_g]v^p+L[v^p]|h|^2_g-4g^{kl}pv^{p-1}\partial_kv\langle
h,\nabla_lh\rangle_g:=(I)+(II)+(III),$$ where:
$$(I)\leq -2v^p|\nabla h|_g^2+c_nv^p|h|_g^4;$$
$$(II)\leq -pv^p|h|_g^4-p(p-1)v^{p-2}|\partial
v|^2_g|h|^2_g+d_npv^{p-1}|\partial v|_g|h|_g^2,$$
$$(III)\leq 4pv^{p-1}|\partial v|_g|h|_g|\nabla h|_g\leq \frac
4{\gamma}\frac p{p-1}v^p|\nabla h|^2_g+\gamma p(p-1)v^{p-2}|\partial
v|^2_g|h|^2_g,$$ for an arbitrary constant $\gamma\in (0,1)$. With
$\eta >0$ to be chosen sufficiently small later, we estimate the
last term in (II):
$$d_npv^{p-1}|\partial v|_g|h|_g^2\leq \eta d_npv^{p-2}|\partial
v|^2_g|h|^2_g+\frac 2{\eta}d_npv^p|h|^2_g.$$ Adding to these
estimates for (I)+(II)+(III) the term $L[B|h|^2_g]$, we have:
$$L[|h|^2_gv^p+B|h|^2_g]\leq [\frac{4p}{\gamma (p-1)}-2-\frac
2{v^p}B]v^p|\nabla h|^2_g+[c_n-p+\frac{c_n}{v^p}B]v^p|h|_g^4$$
$$+[\eta d_np-(1-\gamma)p(p-1)]v^{p-2}|\partial v|^2_g|h|^2_g
+[\frac 2{\eta}d_np]v^p|h|^2_g.$$ Given $\gamma\in (0,1)$ arbitrary,
we choose $p>0$ so large that $(2/\gamma)c_n<p-1$, then $B>0$ so
that:
$$\frac{2p}{\gamma (p-1)}<1+\frac{B}{\bar{v}^p}<1+\frac B{{v}^p}<1+B<\frac{p}{c_n}.$$
In this way we ensure that, in the expression above, the
coefficients in the first two square brackets are negative. Choosing
$\eta>0$ sufficiently small (depending on $\gamma$ and $p$), the
same holds for the third square bracket. Finally, in view of the
second estimate given above for $L[v^p]$, if we add $L[Av^p]$ with
$A>(2/\eta)d_n$ we also take care of the last square bracket (we
also assume $p(p-1)>1/4$), and then:
$$L[Av^p+|h|^2_gv^p+B|h|^2_g]\leq C:=Ad_n^2p^2(\bar{v})^{2(p-1)},$$
concluding the proof. \vspace{.3cm}

\textbf{12. Boundary conditions for the second fundamental form.}
\vspace{.2cm}

Proving global existence for the mean curvature motion of graphs
over time-dependent domains requires estimates for the second
fundamental form.

The simplest form of the evolution equations for $h_{ij}$ and
$H=g^{ij}h_{ij}$ is given in terms of the differential operator on
functions: $L[f]=\partial_tf-tr_gd^2f$.\vspace{.2cm}

The evolution equations for $h_{ij}$ and $H$ are given in Appendix
2. In this section we derive boundary conditions for $h$ and $H$;
the development is similar to work of A. Stahl \cite{Stahl} for MCF
of hypersurfaces intersecting a fixed boundary
orthogonally.\vspace{.2cm}

It is easy to see that $h$ splits on $\partial D(t)$: if $\tau \in
T\partial D(t)$ is a tangential vector field, and $n=n_t$ is the
inner unit normal:
$$h(n,\tau)=\frac 1vd^2w(n,\tau)=\frac 1v(\tau(w_n)-Dw\cdot
\bar{\nabla}_\tau n)=0\mbox{ on }\partial D(t),$$ since $w_n\equiv
\beta_0/\beta$ on the boundary and $\bar{\nabla}_{\tau}n\in
T\partial D(t)$ ($\bar{\nabla}$ is the euclidean connection.) In
particular, it follows that $h(Dw,\tau)=0$ on $\partial D(t)$.
\vspace{.2cm}

\emph{Boundary condition for H.} We derived in section 2 the
equation for the normal velocity of the moving boundary
$\Gamma_t=\partial D(t)$. Letting $\Gamma(\theta,t)$, $\theta\in
S^{n-1}$, be any parametrization of $\Gamma_t$, we find for
$\dot{\Gamma}_n:=\partial_t\Gamma\cdot n$:
$$\dot{\Gamma}_n=-\frac v{w_n}H=-\frac 1{\beta_0}H\quad \mbox{ at
}\partial D(t).$$ Since $\langle N,e_{n+1}\rangle
(t,\Gamma(t))\equiv \beta$ on $\partial D(t)$ we have: $$\langle
\partial_tN,e_{n+1}\rangle=-\langle \partial_kN,e_{n+1}\rangle\dot{\Gamma}^k,$$
where $\partial_kN=-g^{ij}h_{ik}G_j$, with $e_{n+1}$ component:
$$\langle \partial_kN,e_{n+1}\rangle=-g^{ij}w_jh_{ik}=-\frac
1{v^2}h(Dw,\partial_k)=-\frac 1{v^2}w_nh(n,\partial_k).$$ Hence we
find, on $\partial D(t)$:
$$\langle \partial_tN,e_{n+1}\rangle=\frac
{w_n}{v^2}h(n,\dot{\Gamma})=\frac
{w_n}{v^2}\dot{\Gamma}_nh(n,n)=-\beta Hh_{nn}.$$ On the other hand,
using $\partial_tN=-\nabla^{\Sigma}H-Hv^{-1}\nabla^{\Sigma }v$,
combined with the expressions (valid on $\partial D(t)$):
$$\langle \nabla^{\Sigma}H,e_{n+1}\rangle=g^{ij}{H_i}\langle
G_j,e_{n+1}\rangle=g^{ij}H_iw_j=\frac
1{v^2}w_iH_i=\frac{w_n}{v^2}H_n=\beta \beta_0H_n,$$
$$\langle
\nabla^{\Sigma}v,e_{n+1}\rangle=\frac{v_nw_n}{v^2}
=\frac{w_n^2}{v^2}h_{nn}=\beta_0^2h_{nn},$$ we find on $\partial
D(t)$:
$$\langle
\partial_tN,e_{n+1}\rangle=-\beta\beta_0(H_n+\beta_0Hh_{nn}).$$
Comparing these two expressions for $\langle
\partial_tN,e_{n+1}\rangle$ yields:
$$H_n=\frac{\beta^2}{\beta_0}Hh_{nn},$$
a Neumann-type condition for $H$ on $\partial D(t)$.\vspace{.3cm}

\emph{Boundary conditions for $h_{ij}$.} Fix $p\in \partial D(t)$
and let $(\tau_a)$ be an orthonormal frame for $T_p\partial D(t)$
(in the induced metric), satisfying
$\nabla^{\Gamma}_{\tau_a}\tau_b(p)=0$ ($\nabla^{\Gamma}$ is the
connection induced on $\Gamma_t$ by $\bar{\nabla}$, or,
equivalently, by $\nabla$, the Levi-Civita connection of the metric
$g$ in $D(t)$); we extend the $\tau_a$ to a tubular neighborhood so
that $\bar{\nabla}_n\tau_a=0$. Differentiating $h(n,\tau_b)=0$ along
$\tau_a$, we find:
$$(\nabla_{\tau_a}h)(n,\tau_b)=-h(\nabla_{\tau_a}n,\tau_b)-h(n,\nabla_{\tau_a}\tau_b).$$
The second fundamental form ${\cal K}(\tau,\tau')$ of $\Gamma$ in
$(D(t),eucl)$ (equivalently, in $(D(t),g)$) is defined by:
$$\bar{\nabla}_{\tau_a}\tau_b=\nabla^{\Gamma}_{\tau_a}\tau_b+{\cal
K}(\tau_a,\tau_b)n\quad \mbox{ on }\partial D(t).$$ To relate $\cal
K$ to $h_{|\partial D(t)}$, note that since $w=0$ on $\partial
D(t)$:
$$h(\tau_a,\tau_b)=\langle
[\bar{\nabla}_{\tau_a}\tau_b,0],N\rangle=-\bar{\nabla}_{\tau_a}\tau_b\cdot
\frac{Dw}v=-\beta_0{\cal K}(\tau_a,\tau_b).$$ (So we see that
$\Gamma_t$ convex with respect to $n$ corresponds to $\Sigma_t$
\emph{concave} over $D(t)$, as expected). In the appendix we observe
that $\nabla_{\partial_i}\partial_j=(h_{ij}/v)Dw$. Then:
$$\nabla_{\tau_a}\tau_b=\tau_a^i((\tau_b^j)_i\partial_j+\tau_b^j\nabla_{\partial_i}{\partial_j})
=\bar{\nabla}_{\tau_a}\tau_b+\frac 1v\tau_a^i\tau_b^jh_{ij}Dw$$
$$=\nabla^{\Gamma}_{\tau_a}\tau_b+{\cal
K}(\tau_a,\tau_b)n+\frac{w_n}vh(\tau_a,\tau_b)n=(-\frac
1{\beta_0}+\beta_0)h(\tau_a,\tau_b)n=-\frac{\beta^2}{\beta_0}h(\tau_a,\tau_b)n$$
at $p$, given our assumption $\nabla^{\Gamma}_{\tau_a}\tau_b(p)=0$.
We use this immediately to compute, at $p$:
$$\nabla_{\tau_a}n=\langle
\nabla_{\tau_a}n,\tau_b\rangle_g\tau_b=-\langle
n,\nabla_{\tau_a}\tau_b\rangle_g\tau_b=\frac{\beta^2}{\beta_0}|n|^2_gh(\tau_a,\tau_b)\tau_b=\frac
1{\beta_0}h(\tau_a,\tau_b)\tau_b,$$ since
$|n|^2_g=g_{ij}n^in^j=1+w_n^2=\beta^{-2}$ at $p$. We conclude, using
the Codazzi equations:
$$(\nabla_nh)(\tau_a,\tau_b)=(\nabla_{\tau_a}h)(n,\tau_b)=-\frac
1{\beta_0}\sum_ch(\tau_a,\tau_c)h(\tau_c,\tau_b)+\frac{\beta^2}{\beta_0}h(\tau_a,\tau_b)h_{nn}.$$
This can also be written in the form:
$$\beta_0(\nabla_nh)(\tau,\tau')=-(h^{tan})^2(\tau,\tau')+\beta^2h_{nn}h(\tau,\tau').$$
It turns out the expression for covariant derivative of $h$ with
respect to the euclidean connection $\bar{\nabla}$ is exactly the
same (at $\partial D(t)$):
$$\beta_0(\bar{\nabla}_nh)(\tau,\tau')=-(h^{tan})^2(\tau,\tau')+\beta^2
h_{nn}h(\tau,\tau').$$ The reason is that $\nabla_n\tau_a=0$ at the
boundary, also for the $g$-connection:
$$\nabla_n\tau_a=\bar{\nabla}_n\tau_a+n^i\tau_a^j\nabla_{\partial_i}\partial_j=0+\frac
1vh(n,\tau_a)Dw=0,$$ so in fact:
$$(\nabla_nh)(\tau_a,\tau_b)=n(h(\tau_a,\tau_b))=(\bar{\nabla}_nh)(\tau_a,\tau_b).$$

  As done in \cite{Stahl}, we combine this with the
result for $H_n$ to compute $(\nabla_nh)(n,n)$. From:
$$H_n=\nabla_n(tr_gh)=tr_g(\nabla_nh)=\beta^2(\nabla_nh)(n,n)+\sum_a(\nabla_nh)(\tau_a,\tau_a),$$
we find:
$$\beta^2(\nabla_nh)(n,n)=\frac{\beta^2}{\beta_0}Hh_{nn}+\frac
1{\beta_0}|h^{\tan}|^2-\frac{\beta^2}{\beta_0}(H-\beta^2h_{nn})h_{nn}$$
$$=\frac 1{\beta_0}(|h^{tan}|^2+\beta^4h_{nn}^2)=\frac
1{\beta_0}|h|^2_g,$$ since $g^{nn}=\beta^2$ at $\partial D(t)$.
Equivalently:
$$\beta_0(\nabla_nh)(n,n)=\frac 1{\beta^2}|h|^2_g\quad \mbox{ on
}\partial D(t).$$\vspace{.2cm}

It is easy to obtain the corresponding expression for the euclidean
connection. Noting that at $\partial D(t)$:
$$\nabla_nn=\bar{\nabla}_nn+n^in^j\frac 1vh_{ij}Dw=\beta_0h_{nn}n,$$
we find:
$$(\bar{\nabla}_nh)(n,n)=n(h_{nn})=(\nabla_nh)(n,n)+2h(\nabla_nn,n)=(\nabla_nh)(n,n)+2\beta_0h_{nn}^2,$$
so that:
$$\beta_0(\bar{\nabla}_nh)(n,n)=\frac 1{\beta^2}|h|^2_g+2\beta_0^2h_{nn}^2\quad \mbox{ on
}\partial D(t).$$\vspace{.2cm}

It turns out that the expressions just derived, combined with the
maximum principle proved in \cite{Stahl}, are not enough to
establish that concavity is preserved. We derive a suitable maximum
principle in section 13. The result of the next lemma yields a
continuation criterion stated earlier (Prop. 11.4). \vspace{.2cm}

\textbf{Lemma 12.2.} Let $w(y,t)$ be a solution of graph MCM, with
constant-angle boundary conditions, in $E\subset \mathbb{R}^n\times
[0,T)$. Denote by $\cal K$ the second fundamental form of
$\Gamma_t=\partial D(t)$ in $\mathbb{R}^n$. Suppose that, for some
$a_0>0$:
$$\sup\{|A|(y,t)+|\nabla_{\tau}{\cal K}|(y,t);y\in \partial
D(t),t\in [0,T),\tau \in T_y\partial D(t),|\tau|=1\}\leq a_0.$$ Then
also:
$$\sup \{|\nabla A|;y\in \partial D(t), t\in [0,T)\}<\infty.$$

\emph{Proof.} From the boundary conditions computed above for
$\nabla h$, we have at boundary points:
$$|(\nabla_nh)(\tau,\tau)|+|(\nabla_nh)(n,n)|+|(\nabla_{\tau}h)(n,\tau)|\leq
c_0,$$ where $c_0$ depends only on $\beta$ and $a_0$. The remaining
components of $\nabla h$ are:
$$(\nabla_{\tau} h)(n,n)=(\nabla_nh)(\tau,n)\quad \mbox{ and
}(\nabla_{\tau'}h)(\tau,\tau),$$ and since $h^{tan}=-\beta_0{\cal
K}$ at boundary points, the last one is assumed bounded in $[0,T]$.
In addition, on $\partial D(t)$:
$$\tau(H)=\beta^2(\nabla_{\tau}h)(n,n)+2\beta^2h(\nabla_{\tau}n,n)+
\sum_a[(\nabla_{\tau}h)(\tau_a,\tau_a)+2h(\nabla_{\tau}\tau_a,\tau_a)],$$
with all terms on the right bounded, except for the first one. Thus
a bound on $(\nabla_{\tau}h)(n,n)$ would follow from a bound on
$\tau(H)$.But this follows from the uniform gradient estimates (up
to the boundary) of linear parabolic theory, since $H$ is a solution
of (see Appendix 2; $\omega:=Dw/v$):
$$\partial_tH-tr_gd^2H=|h|^2_gH+Hh^2(\omega,\omega)-H^2h(\omega,\omega),\quad {H_n}_{|\partial
D(t)}=\frac{\beta^2}{\beta_0}Hh_{nn},$$ in which all the
coefficients are uniformly bounded in $[0,T)$. The bound depends
only on $a_0$ and the initial data. (The hypotheses of the
proposition imply that the necessary regularity conditions on
$\partial E$ are satisfied.)\vspace{.2cm}

\emph{Finite existence time.}\vspace{.2cm}

In the next section we show that weak concavity at $t=0$ is
preserved by the evolution. Assuming this, it is not difficult to
derive that the flow is defined only for finite time.\vspace{.2cm}

\textbf{Lemma 12.4.} Let $w(y,t),(y,t)\in E\subset {\mathbb
R}^n\times [0,T)$ define a graph MCM $\Sigma_t$ with constant-angle
boundary conditions on a moving boundary. Assume $\Sigma_0$ (and
hence $\Sigma_t$, for all $t$) is weakly concave. Then:

Assume $H_{|t=0}\leq H_0<0$ (where $H_0$ is a negative constant).
Then $T\leq t_*=\frac 1{2H_0^2c_n}$ ( we are assuming $T=\sup \{t\in
[0,T); D(t)\neq \emptyset\}$). Here $c_n>0$ depends only on $n$ and
an upper bound for $v$ in $E$. \vspace{.2cm}

The proof is based on the evolution equation and boundary condition
for $H$ (see Appendix 2: $\omega=Dw/v$):
$$L[H]=|h|^2_gH+Hh^2(\omega,\omega)-H^2h(\omega,\omega),\quad
H_n=(\beta^2/\beta_0)Hh_{nn}.$$
 Since $h^2(\omega,\omega)\geq 0$,  $|h|^2_g \geq
(1/n)H^2$ and (given that $h\leq 0$) $h(\omega,\omega) \geq
|Dw|^2H$, we have:
$$L[H]\leq \frac 1nH^3+|Dw|^2H^3\leq c_nH^3,$$
where $c_n$ depends on $n$ and on $\sup_E |v|$ (already known to be
finite).  Let $\phi(t)$ solve the o.d.e.
$\dot{\phi}=c_n\phi^3,\phi(0)=H_0$:
$$\phi(t)=H_0[1-2c_nH_0^2t]^{-1/2},\quad 0\leq t<t_*:=\frac 1{2H_0^2c_n}.$$
Then with $\psi:=(1/n)(H^2+H\phi+\phi^2)>0$, setting $\chi=H-\phi$:
$$L[\chi]\leq\psi \chi\quad \mbox{ in }E;$$
$$\chi_n=\frac{\beta^2}{\beta_0}(\chi+\phi)h_{nn}\geq
\frac{\beta^2}{\beta_0}\chi\quad \mbox{ on }\partial_lE$$ (since
$\phi<0$ and $h_{nn}\leq 0$). Given that $\chi\leq 0$ at $t=0$, it
follows from the maximum principle that $\chi\leq 0$, or $H\leq
\phi$ in $[0,\min\{T,t_*\})$. This shows $t_*<T$ is impossible,
since $\phi \rightarrow -\infty$ as $t\rightarrow t_*$.
\vspace{.2cm}

\emph{Remark 12.1.} It would be natural to try to show that a
negative upper bound $H_0$ on the mean curvature (at $t=0$) is
preserved, at least under the assumption of concavity.
Unfortunately, the evolution equation for $H$ (under graph m.c.m.)
does not lend itself to a maximum principle argument. Letting
$u:=H-H_0$, we have

$$L[u]= |h|^2_gu+uh^2(\omega,\omega)-u(H+H_0)h(\omega,\omega)+H_0Q\quad \mbox{ in
}E,$$
$$Q:=|h|_g^2+h^2(\omega,\omega)-H_0h(\omega,\omega).$$ At a
point where $u=0$, we would need to show $L[u]\leq 0$. But it is not
true that $Q\geq 0$ at such a point, even when $n=2$. ($u_n\geq 0$
does hold at boundary points.)

\vspace{.4cm}

\textbf{13. A maximum principle for symmetric 2-tensors.}
\vspace{.2cm}

In this section we prove a weak maximum principle for the parabolic
evolution of symmetric two-tensors on bounded euclidean domains,
with moving boundaries and Neumann-type boundary conditions. The
hypotheses are as follows.\vspace{.2cm}

Let $E\subset \mathbb{R}^n\times [0,T]$ be connected, open and
bounded (with $C^2$ boundary), with $D(t)=E\cap ({\mathbb R}^n\times
\{t\})$ bounded, open, connected for each $t\in [0,T]$. Fix $R>0$ so
that $E\subset B_R(0)\times [0,T]$.

On the `lateral boundary' of $E$:
 $$\partial_lE:=\{z=(x,t); t\in [0,T],x\in \partial D(t)\},$$
 we define the \emph{inner} unit normal $n=n_t\in \mathbb R^n$.
 Extend $n_t$ to a vector field in all of $\bar{D}(t)$ (so that it
 is in $C^{2,1}(\bar{E},\mathbb R^n)$, arbitrarily except for the
 requirements that $|n|\leq 1$ pointwise and $\bar{\nabla}_nn=0$ in a tubular neighborhood of
 $\partial D(t)$. (Here $\bar{\nabla}$ denotes the euclidean
 connection, so this requirement can be written $n^i\partial_in^j=0$
 for each $j$.)\vspace{.2cm}

 The assumptions on the coefficients are given next.

 $g=g_t$ is a $t$-dependent Riemannian metric in $\bar{D}(t)$, uniformly equivalent to
 the euclidean metric for $t\in [0,T]$;

 $X=X_t$ is a bounded $t$-dependent vector field in $\bar{D}(t)$,
 satisfying $X\cdot n\geq 0$ for $z\in \partial_lE$;

 $q=q(z,m)$ assigns to each $z\in \bar{E}$ and each $m$ in
 $\mathbb{S}$ (the space of quadratic forms in $\mathbb R^n$) a quadratic form $q\in
 \mathbb{S}$. $q$ is assumed to be $C^{2,1}$ in $z$, locally Lipschitz in
 $m$ (uniformly in $z\in \bar{E}$);

 $b=b(z,m)\in \mathbb{S}$ is defined for $z\in \partial_lE$, with the
 same regularity assumptions.\vspace{.3cm}

 \textbf{Theorem 13.1.} Assume $m\in C^{2,1}(\bar{E};\mathbb{S})$ satisfies in
 $E$ the differential inequality:
 $$\partial_tm_{ij}-tr_gd^2m_{ij}\leq X\cdot dm_{ij}\cdot
 +q_{ij}(\cdot,m(\cdot)),$$
 and on $\partial_lE$ the boundary condition:
 $$n\cdot dm_{ij}(z)\geq b_{ij}(z,m(z)).$$
 Suppose the functions $q$ and $b$ satisfy the following `null
 eigenvector conditions': if, for some $\hat{m}\in \mathbb{S}$,  $V\in {\mathbb R}^n$ is a null
 eigenvector of $\hat{m}$ ($\hat{m}_{ij}V^j=0\forall i$), then, for
 any $z\in \bar{E}$ (resp. any $z\in \partial_lE$):
 $$q_{ij}(z,\hat{m})V^iV^j\leq 0\quad \mbox{ (resp. }b_{ij}(z,\hat{m})V^iV^j\geq
 0).$$
 Then weak concavity of $m$ at $t=0$ is preserved:
 $$m\leq 0\mbox{ in } D(0)\Rightarrow m\leq 0\mbox{ in }\bar{E}.$$
 \vspace{.3cm}

 \emph{Proof.} The assumptions imply there is $K>0$ (depending only on $E$ and
 on the functions $X$, $g$, $n$, $q$ and $b$) satisfying:

 $$|n|_{C^{2,1}(\bar{E})}\leq K,\quad |X(z)|_{eucl}\leq K,\quad |g(z)|+|g^{-1}(z)|\leq K,\quad z\in
 \bar{E};$$
 and if $m,\hat{m}\in C^{2,1}(\bar{E},\mathbb{S})$ satisfy (for some
 $\mu:\bar{E}\rightarrow \mathbb{R}_+$):
 $$-\mu(z)\mathbb{I}\leq m(z)-\hat{m}(z)\leq \mu(z)\mathbb{I}$$
 (where $\mathbb{I}=(\delta_{ij})$ and the inequality of quadratic
 forms has the usual meaning), then also:
 $$q(z,m(z))\leq q(z,\hat{m}(z))+K\mu(z)\mathbb{I},\quad z\in \bar{E},$$
 $$b(z,m(z))\geq b(z,\hat{m}(z))-K\mu(z)\mathbb{I},\quad z\in
 \partial_lE.$$\vspace{.2cm}

 Now define, for $z\in \bar{E}$:
 $$\varphi(z):=-2Kn(z)\cdot x:=2Ks(z),$$
 where we use the euclidean inner product and, on
 $\partial_lE$, $s$ is the `support function' of $\partial D(t)$
 (positive if $D(t)$ is convex and contains the origin). It is clear we may find $M=M(R,K)>0$
 depending only on $K, R$ and $|n|_{C^{2,1}}$ so that:
 $$|\varphi|_{C^{2,1}}\leq M, \quad
 |d\varphi|_g^2+|tr_gd^2\varphi|\leq M,\quad |X\cdot d\varphi|\leq
 M.$$
 We assume also $M\geq K$. Now, given $m$ as in the statement of the theorem and
 given constants $\epsilon>0, \gamma>0$ and
 $\delta >0$, define for $z\in E^{\delta}:=E\cap \{t<\delta\}$:
 $$\hat{m}(z):=m(z)-(\epsilon t+\gamma
 e^{\varphi(z)})\mathbb{I},\quad
 z\in \bar{E}^\delta.$$
 Clearly $\hat{m}\in C^{2,1}(\bar{E}^{\delta}; \mathbb{S})$. We now
 derive the constraints on $\delta$, $\epsilon$ and $\gamma$. It
 will turn out that $\delta$ must be taken small enough (depending only on
 $K,R$), $\epsilon>0$ is arbitrary and $\gamma$ is $\epsilon$ times
 a constant depending only on $K,R$.\vspace{.2cm}

 The following inequalities are easily derived:
 $$q(z,m(z))\leq q(z,\hat{m}(z))+K(\epsilon t+\gamma
 e^{\varphi(z)})\mathbb{I};$$
 $$X\cdot dm=X\cdot d\hat{m}+\gamma (e^{\varphi}X\cdot
 d\varphi)\mathbb{I}\leq X\cdot d\hat{m}+(\gamma
 e^{\varphi}M)\mathbb{I};$$
 $$\partial_t\hat{m}=\partial_tm-\epsilon\mathbb{I}-(\gamma e^{\varphi}\partial_t\varphi)\mathbb{I}
 \leq \partial_tm+(\gamma
 e^{\varphi}M)\mathbb{I}-\epsilon\mathbb{I};$$
 $$tr_gd^2\hat{m}=tr_gd^2\hat{m}-\gamma
 e^{\varphi}(|d\varphi|^2_g+tr_gd^2\varphi)\mathbb{I}
 \geq tr_gd^2m-(\gamma e^{\varphi}M)\mathbb{I}.$$\vspace{.2cm}

 We use this to compute:
 $$\partial_t\hat{m}-tr_gd^2\hat{m}\leq \partial_t
 m-tr_gd^2m+(2\gamma e^{\varphi}M)\mathbb{I}-\epsilon \mathbb{I}$$
 $$\leq q(z,m(z))+X\cdot dm+(2\gamma e^{\varphi}M)\mathbb{I}-\epsilon
 \mathbb{I}$$
 $$\leq q(z,\hat{m}(z))+X\cdot d\hat{m}+K(\epsilon
 t+\gamma e^{\varphi})\mathbb{I}+(3M\gamma e^{\varphi})\mathbb{I}$$
 $$\leq q(z,\hat{m}(z))+X\cdot d\hat{m}+M\epsilon
 t\mathbb{I}+4M\gamma e^{\varphi}\mathbb{I}-\epsilon \mathbb{I}.$$
 We conclude the inequality:
 $$\partial_t\hat{m}-tr_gd^2\hat{m}\leq q(z,\hat{m}(z))+X\cdot
 d\hat{m}-({\epsilon}/2)\mathbb{I}$$ will hold in $E^{\delta}$, provided the constants are selected
 so that, for $z\in E^{\delta}$:
 $$4M\gamma e^{\varphi(z)}+M\epsilon t\leq {\epsilon}/2.\qquad
 (A)$$\vspace{.2cm}
 Turning to boundary points $z=(x,t)\in \partial_lE$, note that
 $d_n\varphi=-2K$, so that:
 $$d_n\hat{m}(z)=d_nm(z)-(\gamma
 e^{\varphi(z)}d_n\varphi(z))\mathbb{I}\geq b(z,m(z))-(\gamma
 e^{\varphi(z)}d_n\varphi(z))\mathbb{I}$$
 $$\geq b(z,\hat{m}(z))-K(\epsilon t+\gamma
 e^{\varphi(z)})\mathbb{I}-(\gamma
 e^{\varphi(z)}d_n\varphi(z))\mathbb{I}$$
 $$\geq b(z,\hat{m}(z))+K(\gamma e^{\varphi(z)}-\epsilon
 t)\mathbb{I},$$
 so that the inequality:
 $$d_n\hat{m}(z)\geq b(z,\hat{m}),\quad z\in \partial_lE^{\delta},$$
 will hold provided the constants are chosen so that, on
 $\partial_lE^{\delta}$:
 $$\epsilon t\leq \gamma e^{\varphi(z)}.\qquad (B)$$
 Bearing in mind that, on $E$: $e^{-2KR}\leq e^{\varphi(z)}\leq
 e^{2KR}$, it is not hard to arrange for (A) and (B) to hold, or
 equivalently, for:
 $$\epsilon t\leq \gamma e^{\varphi(z)},\qquad 10M\gamma
 e^{\varphi(z)}\leq \epsilon.$$
 Given $\epsilon >0$, define $\gamma$ so that $10M\gamma
 e^{2KR}=\epsilon$. Then the second inequality holds, and so will
 the first, provided:
 $$\epsilon t\leq \gamma e^{-2KR}=(\epsilon/10M)e^{-4KR},$$
 which is true for any $\epsilon>0$, if $\delta$ is defined via
 $\delta:=e^{-4KR}/10M$ (recall $t\in [0,\delta]$).\vspace{.2cm}

  Note that, since $m\geq 0$ at $t=0$, it follows that $\hat{m}$ is negative-definite
  at $t=0$, and hence also
 for small time, and we \emph{claim} that this persists
 throughout $\bar{E}^{\delta}$, so that (letting $\epsilon
 \rightarrow 0$) $m\leq 0$ in $\bar{E}^{\delta}$. Restarting the
 argument at $t=\delta$, we see this is enough to prove the theorem.
 \vspace{.2cm}

 To prove this claim, suppose (by contradiction) $\hat{m}$ acquires
 a null eigenvector $0\neq V\in {\mathbb R}^n$ at a point
 $z_1=(x_1,t_1)\in \bar{E}^{\delta}$, with $t_1\in (0,\delta]$ the
 first time this happens.

 Let $\hat{f}(z):=\hat{m}_{ij}V^iV^j, \quad z\in E^{\delta}$ (that
 is, we `extend' $V$ to $E^{\delta}$ as a constant vector.) It
 follows from the preceding that $\hat{f}$ satisfies in
 $E^{\delta}$:
 $$\partial_t\hat{f}\leq tr_gd^2\hat{m}_{ij}V^iV^j+X\cdot
 d\hat{m}_{ij}V^iV^j+q_{ij}(\cdot,\hat{m})V^iV^j-\frac{\epsilon}2|V|^2_{eucl}.$$
 Noting that $tr_gd^2\hat{m}_{ij}V^iV^j=tr_gd^2\hat{f}$ and $X\cdot
 \hat{m}_{ij}V^iV^j=d\hat{f}\cdot X$, and using the null eigenvector
 condition for $q$, we find that $\hat{f}$ satisfies in $E^{\delta}$
 the strict inequality:
 $$\partial_t\hat{f}<tr_gd^2\hat{f}+d\hat{f}\cdot X.$$
 This shows $x_1$ cannot be an interior point of $D(t_1)$, for then
 (as a first-time interior maximum point for $\hat{f}$) we would
 have $tr_gd^2\hat{f}(z_1)\leq 0$ and $d\hat{f}(z_1)=0$,
 contradicting $\partial_t\hat{f}(z_1)\geq 0$. Thus $x_1\in \partial
 D(t_1)$. Since $\hat{f}$ satisfies the differential inequality just
 stated and $z_1=(x_1,t_1)$ is a first-time boundary maximum in
 $\bar{E}^{\delta}$, the parabolic Hopf lemma implies
 $d_n\hat{f}(z_1)<0$. On the other hand, as seen above:
 $$d_n\hat{f}= d_n\hat{m}_{ij}V^iV^j\geq
 b_{ij}(z_1,\hat{m}(z_1))V^iV^j\geq 0,$$
 from the boundary null-eigenvector condition. This contradiction
 concludes the proof.\vspace{.3cm}

 \textbf{Corollary 13.2.} Suppose $m\in C^{2,1}(\bar{E},\mathbb{S})$
 satisfies the same differential inequality, with the same
 hypotheses on the coefficients as in the theorem (including the
 null eigenvector condition for $q$), and the boundary conditions:
 $$m(z)(n,\tau)=0,\quad \forall z=(x,t)\in \partial_lE,\tau\in
 T_x\partial D(t);$$
 $$n^in^jd_n m_{ij}(z)=(\bar{\nabla}_nm)(n,n)\geq b_{nn}(z,m(z));$$
 $$\tau^i\tau^jd_nm_{ij}=(\bar{\nabla}_nm)(\tau,\tau)\geq
 b^{tan}(z,m(z))(\tau,\tau),\quad \tau \in T_x\partial D(t),$$
 for functions $b_{nn}(z,\hat{m})$ from $E\times \mathbb{S}$ to $\mathbb{R}$ and $b^{tan}$ assigning
 to $(z,\hat{m})$, $z=(x,t)$, a quadratic form in $T_{x}\partial D(t)$.
 Suppose $b_{nn}\geq 0$ in $E\times \mathbb{S}$ and $b^{tan}$ satisfies:
 $$\hat{m}(\tau,\tau)=0\mbox{ for some }\tau\in T_x\partial
 D(t)\Rightarrow b^{tan}(z,\hat{m})(\tau,\tau)\geq 0.$$
 Then, as in the theorem, concavity is preserved:
 $$m\leq 0\mbox{ at }t=0\Rightarrow m\leq 0\mbox{ in
 }\bar{E}.$$\vspace{.2cm}

 \emph{Proof.} This is proved as the theorem, with the following
 change in the last part of the proof: if $0\neq V\in \mathbb{R}^n$
 is a null eigenvector of $\hat{m}$ (defined as in the proof of the
 theorem) at a boundary point $z_1=(a_1,t_1)\in \partial_l E$,
 write:
 $$V=V^nn+V^T,\quad V^T\in T_{x_1}\partial D(t_1).$$
 Assume first $V^n\neq 0$. Then (noting that $\hat{m}$ splits at the
 boundary if $m$ does), we see that $n$ is a null eigenvector of
 $\hat{m}$ at $z_1$, so we define
 $\hat{f}(z)=m_{ij}(\hat{z})n^i(z_1)n^j(z_1)$ and repeat the
 argument. At $z_1$,
 $(\bar{\nabla}_n\hat{m})(n,n)=b_{nn}(z_1,\hat{m}(z_1))\geq 0$ leads to a
 contradiction with the parabolic Hopf lemma, as before.

 If $V^n=0$, then $V^T\in T_{x_1}\partial D(t_1)$ must be a null eigenvector of $\hat{m}$ at
 the boundary point $z_1$, and then we run the argument with
 $\hat{f}(z)=\hat{m}(z)(V^T,V^T)$, leading to a contradiction,
 as before.\vspace{.3cm}

 \textbf{Corollary 13.3.} For MCM of graphs with constant-angle
 boundary conditions, weak concavity is preserved:
 $$h\leq 0\mbox{ at }t=0\Rightarrow h\leq 0\mbox{ in }\bar{E}.$$

 \emph{Proof.} The conditions of the theorem hold, and the
 expressions obtained for $\bar{\nabla}_nh$ in the preceding section easily imply that the
 boundary conditions in Corollary 13.2 are satisfied; hence the
 claim follows from Corollary 13.2.\vspace{.3cm}

 \emph{Remark 13.1.} It seems plausible that a slightly different version
 of the result in this section could be used to strengthen the
 conclusions of \cite{Stahl}. This is currently being considered.

\vspace{.5cm}

\textbf{14. An improved continuation criterion.}\vspace{.2cm}

In this section we improve the continuation criterion: if
$\sup_E|h|_g=a_0$ is finite, the solution can be continued past $T$.

The argument given below works in all dimensions, but for simplicity
of notation we deal here only with the two-dimensional case:
$\Sigma_t$ is a surface, the moving boundary $\Gamma_t$ is a curve
in $\mathbb R^2$. Assuming such a bound on $|h|_g$, given the
results in section 12 all we have to do is bound
$(\nabla_{\tau}h)(\tau,\tau)$ and $(\nabla_{\tau}h)(n,n)$, where
$\tau=\tau_t$ is a unit vector field tangent to $\Gamma_t$. At the
moving boundary: $H=\beta^2 h(n,n)+h(\tau,\tau)$, and we already
showed $\tau(H)$ is bounded, so it suffices to bound one of these
quantities.\vspace{.2cm}

We adopt the \emph{notation}: $f\sim g$ if $f-g$ is bounded in $E$
by constants depending only on the initial data and
$a_0$.\vspace{.2cm}

Consider the vector fields in $D(t)\subset \mathbb{R}^2$:
$$\omega=\frac 1v[w_1,w_2], \quad
\tilde{\omega}=v\omega^{\perp}=[-w_2,w_1].$$ It is easy to verify
the following:
$$\langle \omega,\tilde{\omega}\rangle_g=0,\quad
|\omega|^2_g=|\tilde{\omega}|^2_g=|Dw|^2_e:=w_1^2+w_2^2.$$ Thus we
may think of $\{\omega, \tilde{\omega}\}$ as a `conformal
pseudo-frame' ($\omega$ and $\tilde{\omega}$ vanish when $Dw=0$),
defined on all of $D(t)$. Moreover, at the boundary $\partial D(t)$:
$$\omega=\beta_0n,\quad
\tilde{\omega}=\frac{\beta_0}{\beta}n^{\perp}:=\frac{\beta_0}{\beta}\tau,$$
where $\{\tau,n\}$ is an euclidean-orthonormal frame along
$\Gamma_t$. Thus $\omega$ and $\tilde{\omega}$ supply `canonical'
extensions of $n,\tau$ to the interior of $D(t)$, as uniformly
bounded vector fields.\vspace{.2cm}

Recall the boundary conditions for $h$:
$$h(\omega,\tilde{\omega})=0,\quad (\nabla_nh)(\omega,\omega)\sim
0,\quad (\nabla_nh)(\tilde{\omega},\tilde{\omega})\sim 0\quad \mbox{
on }\partial D(t).$$

These give the boundary conditions for the components
$h_{11},h_{12},h_{22}$ of $h$ in the standard basis of
$\mathbb{R}^2$. As shown in appendix 2, these three functions are
solutions of a linear parabolic system in $E\subset {\mathbb
R}^2\times [0,T]$ (a non-cylindrical domain), with bounded
coefficients:
$$\partial_t(h_{ij})-g^{kl}(h_{ij})_{kl}+2h_i^kd_{\omega}(h_{kj})+2h_j^kd_{\omega}(h_{ik})=C_{ij}.$$

The boundary conditions also have bounded coefficients:
$$n^1n^2(h_{22}-h_{11})+[(n^1)^2-(n^2)^2]h_{12}=0,$$
$$(n^1)^2d_n(h_{11})+2n^1n^2d_n(h_{12})+(n^2)^2d_n(h_{22})=b_1,$$
$$(n^2)^2d_n(h_{11})-2n^1n^2d_n(h_{12})+(n^1)^2d_n(h_{22})=b_2.$$

The only thing left to do is to argue that this set of linear
equations and boundary conditions define a parabolic system. Then it
follows from the `global gradient bounds' of linear theory that also
the tangential derivatives $d_{\tau}h_{ij}$ are bounded on
$\partial_lE$, which leads quickly to the desired
conclusion.\vspace{.2cm}

We need to verify the `complementarity conditions' hold for this
system, so we proceed as in Section 6 (up to a point.) Fix a point
$z_0=(y_0,t_0)\in \partial_lE$ and `freeze coefficients' there.
Consider a manifold-with-boundary chart $(y,t)\mapsto
(\rho,\sigma,s)$ mapping a neighborhood of $z_0$ in $E$ to
$\{\rho>0\}\times {\mathbb R}\times {\mathbb R}_+$. Here $\rho$ is
the coordinate normal to $\partial D(t)$, $\sigma$ parametrizes
$\partial D(t)$ and slices $\{s=const.\}$ correspond to $\{
t=const.\}$.

Let $\tilde{h}_{ij}(\rho,\sigma,s)=h_{ij}(y,t)$ be the unknown
functions in the new coordinates. The corresponding system is:
$$\partial_s\tilde{h}_{ij}-(\beta^2 (\tilde{h}_{ij})_{\rho
\rho}+(\tilde{h}_{ij})_{\sigma
\sigma})+c^k_i(\tilde{h}_{jk})_{\rho}+c^k_j(\tilde{h}_{ik})_{\rho}=\varphi_{ij}.$$
Here the $c_i^k$ are constants. The boundary conditions can also be
easily written down (freeze the $n^i$ to their value $n^i_0$ at
$z_0$ and replace $d_n(h_{ij})$ by $\tilde{h}_{ \rho}$.) It is
natural to consider the linear transformation of the unknown
functions:
$$(\tilde{h}_{11},\tilde{h}_{12},\tilde{h}_{22})\mapsto
(f_{11},f_{12},f_{22}),\quad f_{ij}=f_{ij}(\rho,\sigma,s)$$ defined
by:
$$f_{11}=(n_0^1)^2\tilde{h}_{11}+2n_0^1n_0^2\tilde{h}_{12}+(n_0^2)^2\tilde{h}_{22}$$
$$f_{22}=(n_0^2)^2\tilde{h}_{11}-2n_0^1n_0^2\tilde{h}_{12}+(n_0^1)^2\tilde{h}_{22}$$
$$f_{12}=(n_0^1n_0^2)(\tilde{h}_{22}-\tilde{h}_{11})+[(n_0^1)^2-(n_0^2)^2]\tilde{h}_{22}.$$
Since the principal part of the linear system for the
$\tilde{h}_{ij}$ is diagonal, the principal part of the system for
$f_{ij}$ is exactly the same (while the lower-order terms have
different values):
$$\partial_sf_{ij}-(\beta^2 (f_{ij})_{\rho
\rho}+(f_{ij})_{\sigma
\sigma})+\bar{c}^k_i(f_{jk})_{\rho}+\bar{c}^k_j(f_{ik})_{\rho}=\bar{\varphi}_{ij},$$
This linear transformation is invertible, with inverse given by:
$$\tilde{h}_{11}=(n_0^2)^2f_{22}+2n_0^1n_0^2f_{12}+(n_0^1)^2f_{11}$$
$$\tilde{h}_{22}=(n_0^1)^2f_{22}-2n_0^1n_0^2f_{12}+(n_0^2)^2f_{11}$$
$$\tilde{h}_{12}=n_0^1n_0^2(f_{11}-f_{22})+[(n_0^2)^2-(n_0^1)^2]f_{12}$$

Thus the original system frozen at $z_0$ (for the $\tilde{h}_{ij}$),
with its boundary conditions, satisfies the complementarity
condition if an only if the same holds for the $f_{ij}$ system, with
the transformed boundary conditions. But these take a very simple
form:
$${f_{12}}_{|\rho=0}=0,\quad
d_{\rho}(f_{11})_{|\rho=0}=\bar{b}_1,\quad
d_{\rho}(f_{22})_{|\rho=0}=\bar{b}_2.$$ These are standard Dirichlet
(resp. Neumann) boundary conditions for a standard $3\times 3$
parabolic system (decoupled to highest order). Hence the original
system (for the $h_{ij}$) with boundary conditions satisfies
`complementarity' at each point of $\partial_l E$.\vspace{.2cm}

In particular, the global gradient estimates hold for the linear
system (with uniformly bounded coefficients and boundary conditions)
in the unknowns $h_{ij}$, and we have:
$$|d_{\tau}h_{ij}|\leq M\quad \mbox{ in }\partial_lE,$$
for any tangential unit vector field $\tau$, for some $M$ depending
only on $a_0$ and the initial data. This clearly implies bounds on
$(\nabla_{\tau}h)(\tau,\tau)$ and $(\nabla_\tau h)(n,n)$. Combining
with the results in section 12, we have the following conclusion:
\vspace{.2cm}

\textbf{Proposition 14.1.} Assume the maximal existence time
$T_{max}$ is finite. Then:
$$\limsup_{t\rightarrow T_{max}}\sup_{\partial
D(t)}|h|_g=\infty.$$\vspace{.2cm}

\emph{Remark.} It should be clear that the argument works in all
dimensions; this will be included in the final version of the
paper.\vspace{.5cm}

\textbf{15. Final comments.}\vspace{.2cm}

1. The main step missing for the global existence result
$$\lim_{t\rightarrow T_{max}}diam(\Sigma_t)=0$$ (in the concave case)
is showing that a lower bound on diameter gives an upper bound for
$|h|_g$. This may follow from properties of the support function
(based on a point in $\mathbb{R}^n$ common to all $\Sigma_t$), but
remains to be addressed (work in progress). If confirmed, this would
correspond to Theorem 1 in \cite{LensSeminar} for lens-type curve
networks. An issue apparently completely unexplored in dimensions
above 1 is existence-uniqueness of `homothetic solutions' for this
problem.\vspace{.2cm}

2. We state here the local existence theorem for configurations of
graphs over domains with moving boundaries. In this setting, a
\emph{triple junction configuration} consists of three embedded
hypersurfaces $\Sigma^1,\Sigma^2,\Sigma^3$ in $\mathbb{R}^{n+1}$,
graphs of functions $w^I$ defined over time-dependent domains
$D^1(t), D^2(t)\subset \mathbb{R}^n$ ($D^1$ covered by one graph,
$D^2$ by two graphs), satisfying the following conditions: (1) The
$\Sigma^I$ intersect along an $(n-1)$-dimensional graph $\Lambda(t)$
(the `junction'), along which the upward unit normals satisfy the
relation: $N_1+N_2=N_3$. (2) If a fixed support hypersurface
$S\subset {\mathbb R}^{n+1}$ is given (also a graph, not necessarily
connected), the $\Sigma^I$ intersect $S$ orthogonally.

Topologically, in the case of bounded domains one has the following
examples: (i) (`lens' type) 2 disks (or two annuli) covering
$D^2(t)$ and one annulus covering $D^1(t)$; (ii) (`exterior' type)
two annuli covering $D^2(t)$ and one disk covering $D^1(t)$. The
boundary component of the annuli disjoint from the junction
intersects the support hypersurface $S$ orthogonally for each
$t$.\vspace{.2cm}

 Let $\Sigma_0^I$ ($I=1,2,3$) be graphs of $C^{3+\alpha}$
functions over $C^{3+\alpha}$ domains $D_0^1,D_0^2\subset
\mathbb{R}^n$, defining a triple junction configuration and
satisfying the compatibility condition for the mean curvatures on
the common boundary $\Gamma_0$ of $D^1_0$ and $D^2_0$:
$$H^1+H^2=H^3.$$
Then there exists $T>0$ depending only on the initial data, and
functions $w^I\in C^{2+\alpha,1+\alpha/2}(Q^I)$, $Q^I\subset
{\mathbb R}^n\times [0,T)$, so that the graphs of
$w^I(.,t):D^I(t)\rightarrow {\mathbb R}$ define a triple junction
configuration for each $t\in [0,T)$, moving by mean
curvature.

The proof will be given elsewhere. \vspace{.2cm}

3. An interesting issue we have not addressed here is whether one
has breakdown of uniqueness for initial data of lower regularity, or
if the `orthogonality condition' at the junction is removed. For
curve networks, non-uniqueness has been considered in
\cite{MazzeoSaez}; but neither a drop in regularity (from initial
data to solution, in H\"{o}lder spaces) nor the orthogonality
condition play a role in the case of curves.\vspace{.5cm}

\emph{Appendix 1}: \textbf{Proof of lemma 4.1.}\vspace{.2cm}

Throughout the proof, $n$ denotes the inner unit normal at $\partial
D$, extended to a tubular neighborhood $\cal{N}$ so that $D_nn=0$.
Since $D$ is  uniformly $C^{3+\alpha}$, if follows that $n\in
C^{2+\alpha}(\partial D)$, with uniform bounds. Denote by $\rho$ the
distance to the boundary (so $D\rho=n$ in $\cal{N}$). Let $\zeta\in
C^3(\bar{D})$ be a cutoff function, with $\zeta\equiv 1$ in ${\cal
N}_1\subset {\cal N}$, $\zeta \equiv 0$ in $D\setminus {\cal
N}$.\vspace{.2cm}

We find $\varphi$ of the form:
$$\varphi(x)=x+\zeta(x)f(x)n(x)$$
with $f\in C^{2+\alpha}({\cal N})$. The 1-jet conditions on
$\varphi$ at $\partial D$ translate to the conditions on $f$:
$$f_{|\partial D}=0,\quad Df_{|\partial D}=0,\quad
D^2f(n,n)_{|\partial D}={\Delta f}_{|\partial D}=h.$$ Now
use:\vspace{.2cm}

\emph{Lemma A.1.} Let $D$ be a uniformly $C^{3+\alpha}$ domain with
boundary distance function $\rho>0$. Let $h\in C^{\alpha}(\partial
D)$ be a bounded function. Then there exists an extension $g\in
C^{\infty}(D)\cap C(\bar{D})$ so that $g_{|\partial D}=h$,
$\sup_{\bar{D}}|g|\leq \sup_{\partial D}|h|$ and $\rho^2 g\in
C^{2+\alpha}(\bar{D})$.

Given this lemma, all we have to do is set $f=(1/2)\rho^2g$, which
clearly satisfies all the requirements (in particular, $\Delta f=h$
at $\partial D$.)\vspace{.2cm}

To verify that $\varphi$ is a diffeomorphism, it suffices to check
that $|\zeta fn|_{C^1}$ (in ${\cal N}\subset \{ \rho<\rho_0\}$) is
small if $\rho_0$ is small. This is easily seen:

$$|\zeta f n|_{C^0}\leq (1/2)\rho_0^2|g|_{C^0};$$
$$|D\zeta|\leq c\rho_0^{-1}\Rightarrow |fD\zeta|\leq c\rho_0
|g|_{C^0}.$$
$$|Df|\leq (1/2)\rho_0^{\alpha}||g||_{C^{2+\alpha}(\bar{D})}$$
on $\cal{N}$, since $Df\in C^{1+\alpha}(\bar{D})$ and $Df_{|\partial
D}=0$. And finally, with ${\cal A}$ the second fundamental form of
$\partial D$:
$$|Dn|\leq |{\cal A}|_{C^0}\Rightarrow
|fDn|\leq(1/2)\rho_0^2|g|_{C^0}|{\cal A}|_{C^0}.$$\vspace{.2cm}

A word about Lemma A.1. (This is probably in the literature, but I
don't know a reference.) If $D$ is the upper half-space, we solve
$\Delta g=0$ in $D$ with boundary values $h$. Then the estimate
$$[D^2(\rho^2 P*h)]^{(\alpha)}(\bar{D})\leq c
|h|_{C^{\alpha}(\partial D)}$$ follows by direct computation with
the Poisson kernel $P$; for the rest of the norm, use interpolation.
Then transfer the estimate to a general domain using `adapted local
charts', in which $\rho$ in $D$ corresponds to the vertical
coordinate in the upper half-space. (It is easy to see that at each
boundary point there is a $C^{2+\alpha}$ adapted chart, with uniform
bounds.)\vspace{.5cm}

\noindent \emph{\textbf{Appendix 2: Evolution equations for the
second fundamental form.}}\vspace{.2cm}

We consider mean curvature motion of graphs:
$$G(y,t)=[y,w(y,t)],\quad y\in D(t)\subset \mathbb{R}^n,$$
$$w_t=g^{ij}w_{ij}=vH,\quad v=\sqrt{1+|Dw|^2}.$$
In this appendix we include evolution equations for geometric
quantities, in terms of the operators:
$$\partial_t-\Delta_g,\quad \quad L=\partial_t-tr_gd^2.$$
It is often convenient to use the vector field in $D(t)$:
$$\omega:=\frac 1v Dw.$$
Since $-\omega$ is the ${\mathbb R^n}$ component of the unit normal
$N$ and $L[N]=|h|^2_gN$, we have:
$$L[\omega^i]=|h|^2_g\omega^i, \quad |h|^2_g:=g^{ik}g^{jl}h_{ij}h_{kl}.$$
Here  $h=(h_{ij})$ is the pullback to $D(t)$ of the second
fundamental form $A$:
$$h(\partial_i,\partial_j)=h_{ij}=A(G_i,G_j)=\frac
1vw_{ij}.$$

First, denoting by $\nabla$ the pullback to $D(t)$ of the induced
connection $\nabla^{\Sigma}$ (that is,
$G_*(\nabla_XY)=\nabla_{G_*X}^{\Sigma}G_*Y$ for any vector fields
$X,Y$ in $D(t)$), and using the definition:
$$\nabla^{\Sigma}_{G_i}G_j=G_{ij}-\langle G_{ij},N\rangle
N=[0,w_{ij}]-\frac
1{v^2}w_{ij}[-Dw,1]=\frac{w_{ij}}{v^2}[Dw,|Dw|^2]=\frac{w_{ij}}{v^2}G_*Dw,$$
we conclude:
$$\nabla_{\partial_i}\partial_j=\frac 1vh_{ij}Dw=h_{ij}\omega.$$

From this one derives easily a useful expression relating the
Laplace-Beltrami operator and the operator $tr_gd^2$ acting on
functions:
$$\Delta_gf=tr_gd^2f-\frac Hvw_mf_m=tr_gd^2f-Hd_{\omega}f.$$
\vspace{.2cm}

We also have, for the covariant derivatives of $h$ with respect to
the euclidean connection and to $\nabla=\nabla^g$:
$$\partial_m(h_{ij})=\nabla_mh_{ij}+[h_{jm}h_{ik}+h_{im}h_{jk}]\omega^k.$$
(Here $\nabla h$ is the symmetric $(3,0)$-tensor with components:
$\nabla_mh_{ij}=(\nabla_{\partial_m}h)(\partial_i,\partial_j)$.)\vspace{.2cm}

Iterating this and taking $g$-traces yields (using the Codazzi
identity and the easily verified relation
$\partial_i\omega^k=h^k_i:=g^{jk}h_{ij}$):
$$tr_gd^2(h_{ij})=g^{mk}\partial_m(\partial_k(h_{ij}))=g^{mk}(\nabla^2_{\partial_m,\partial_k}h)(\partial_i,\partial_j)$$
$$+H\nabla_{\omega}h_{ij}+2[h_i^k\nabla_kh_{jp}+h_j^k\nabla_kh_{ip}]\omega^p+[H_ih_{jp}+H_jh_{ip}]\omega^p$$
$$+2[h_{ip}(h^2)_{jq}+(h^2)_{ip}h_{jq}+Hh_{ip}h_{jq}]\omega^p\omega^q+2(h^3)_{ij}+2(h^2)_{ij}h(\omega,\omega).$$

Here the powers $h^2$ and $h^3$ of $h$ are the symmetric 2-tensors
defined used the metric:
$$(h^2)_{ij}:=g^{kp}h_{ik}h_{pj}=h^k_ih_{pj},\quad (h^3)_{ij}:=g^{kp}g^{lq}h_{ik}h_{pl}h_{qj}.$$
Note also that:
$$[h_i^k\nabla_kh_{jp}+h_j^k\nabla_kh_{ip}]\omega^p=\nabla_{\omega}(h^2)_{ij},$$
using the Codazzi identity.\vspace{.2cm}

\emph{Evolution equations for $h$.}

Starting from $G_t=vHe_{n+1}=H(N+\frac 1v[Dw,|Dw|^2])=HN+HG_*\omega$
and $N_t=-\nabla^{\Sigma} H-Hv^{-1}\nabla^{\Sigma}v$ (where
$\nabla^{\Sigma}f=g^{ij}f_jG_i$ and $\nabla f=g^{ij}f_j\partial_i$)
we have:
$$\partial_t(h_{ij})=\langle (HN)_{ij},N\rangle-\langle
G_{ij},\nabla^{\Sigma}H\rangle-\frac Hv\langle G_{ij},\nabla^\Sigma
v\rangle+\langle (HG_*\omega)_{ij},N\rangle.$$ Using the easily
derived facts:
$$\langle N_{ij},N\rangle=-h^2(\partial_i,\partial_j),$$
$$H_{ij}-\langle G_{ij},\nabla^{\Sigma}H\rangle=(\nabla
dH)(\partial_i ,\partial_j),$$
$$\frac 1v\langle G_{ij},\nabla^{\Sigma}
v\rangle=h(\omega,\omega)h_{ij},$$ we obtain:
$$\partial_t(h_{ij})=(\nabla
dH)(\partial_i,\partial_j)-Hh^2(\partial_i,\partial_j)-Hh(\omega,\omega)h_{ij}+
\langle (HG_*\omega)_{ij},N\rangle,$$ where:
$$\langle (HG_*\omega)_{ij},N\rangle=H_i\langle
(G_*\omega)_j,N\rangle+H_j\langle (G_*\omega)_i,N\rangle+H\langle
(G_*\omega)_{ij},N\rangle.$$
 To identify the terms, computation shows that:
$$\langle (G_*\omega)_i,N\rangle=h(\omega,\partial_i),$$
and hence, using also:
$$\nabla^{\Sigma}_{G_i}(G_*\omega)=G_*(\nabla_{\partial_i}\omega),\quad
\nabla_{\partial_i}\omega=(h_i^p+\omega^qh_{iq}\omega^p)\partial_p=\sum_ph_{ip}\partial_p,$$
we obtain (using
$\omega^k\partial_j(h_{ik})=\nabla_{\omega}h_{ij}+2h(\partial_i,\omega)h(\partial_j,\omega)$):
$$\langle
(G_*\omega)_{ij},N\rangle=\partial_j(\omega^kh_{ik})-\langle
\nabla_{G_i}^{\Sigma}(G_*\omega),\partial_jN\rangle
=h_j^kh_{ik}+\omega^k\partial_j(h_{ik})+h(\partial_j,\nabla_{\partial_
i}\omega)$$
$$=(\nabla_{\omega}h)_{ij}+(h^2)_{ij}+2h(\omega,\partial_i)h(\omega,\partial_j)+\sum_ph_{ip}h_{jp}$$
$$=(\nabla_{\omega}h)_{ij}+2(h^2)_{ij}+3h(\omega,\partial_i)h(\omega,\partial_j),$$
since
$\sum_ph_{ip}h_{jp}=(h^2)_{ij}+h(\omega,\partial_i)h(\omega,\partial_j).$
 Combining all the terms yields the result:
$$\partial_t(h_{ij})=(\nabla
dH)(\partial_i,\partial_j)+H\nabla_{\omega}h_{ij}+H_ih(\omega,\partial_j)+H_jh(\omega,\partial_i)$$
$$+H(h^2)_{ij}+3Hh(\omega,\partial_i)h(\omega,\partial_j)-Hh(\omega,\omega)h_{ij}.$$\vspace{.2cm}

From this expression and Simons' identity (in tensorial form):
$$\nabla dH=\Delta_gh+|h|^2_gh-Hh^2,$$
we obtain easily a tensorial `heat equation' for $h$:
$$[(\partial_t-\Delta_g)h]_{ij}=H\nabla_{\omega}h_{ij}+H_ih(\omega,\partial_j)+H_jh(\omega,\partial_i)$$
$$+|h|^2_gh_{ij}+3Hh(\partial_i,\omega)h(\partial_j,\omega)-Hh(\omega,\omega)h_{ij}.$$\vspace{.2cm}

Using the earlier computation relating $\Delta_gh$ (the tensorial
Laplacian  of $h$) and $tr_gd^2h$, we obtain from this the evolution
equation in terms of $L$:
$$L[h_{ij}]=-2\nabla_{\omega}(h^2)_{ij}+C_{ij},$$
$$C_{ij}:=-2[h(\partial_i,\omega)h^2(\partial_j,\omega)+h^2(\partial_i,\omega)h(\partial_j,\omega)]-2(h^3)_{ij}-2(h^2)_{ij}h(\omega,\omega)$$
$$+|h|^2_gh_{ij}+Hh(\partial_i,\omega)h(\partial_j,\omega)-Hh(\omega,\omega)h_{ij}.$$\vspace{.2cm}

\emph{Time derivatives and evolution equations for $\omega$ and
$g$.}

It is sometimes convenient to use the `Weingarten operator':
$$S(X):=S(X^i\partial_i)=h^i_jX^j\partial_i.$$

The time derivative of $\omega$ is simply minus the time derivative
of the $\mathbb{R}^n$ component of $N$. In addition, one computes
easily that $\frac{\nabla v}v=S(\omega)$, so we have:
$$\partial_t\omega=\nabla H+\frac Hv \nabla v=\nabla H+HS(\omega).$$

For the metric and `inverse metric' tensors we have: from
$\partial_tg_{ij}=(w_iw_j)_t$ and $w_{it}=(vH)_i$:
$$\partial_tg_{ij}=v^2(H_i\omega^j+H_j\omega^i)+v^2H(h(\omega,\partial_i)\omega^j+h(\omega,\partial_j)\omega^i),$$
and then, using $\partial_tg^{ij}=-g^{ik}\partial_tg_{kl}g^{lj}$:
$$\partial_tg^{ij}=- [(\nabla H)^i\omega^j+(\nabla
H^j)\omega^i]-H[S(\omega)^i\omega^j+S(\omega)^j\omega^i].$$

Since we know the evolution equation of $\omega$, it is easy to
obtain that of $g^{ij}$:
$$L[g^{ij}]=-L[\omega^i\omega^j]=-L[\omega^i]\omega^j+2g^{kl}(\partial_k\omega^i)(\partial_l\omega^j)-\omega^iL[\omega^j].$$
Using $\partial_k\omega^i=h_k^i$, we find:
$$L[g^{ij}]=-2|h|^2_g\omega^i\omega^j+2(h^2)^{ij}.$$
It is also easy to see that
$\partial_kg^{ij}=-(h_k^i\omega^j+h_k^j\omega^i)$.\vspace{.2cm}

\emph{Evolution of mean curvature.}

To compute the evolution equation for $H=g^{ij}h_{ij}$, we just need
to remember $g^{ij}$ is time-dependent:
$$(\partial_t-\Delta_g)H=(\partial_tg^{ij})(h_{ij})+tr_g[(\partial_t-\Delta_g)h]
=-2h(\nabla
H,\omega)-2Hh^2(\omega,\omega)+tr_g[(\partial_t-\Delta_g)h].$$ The
result is:
$$(\partial_t-\Delta_g)H=Hd_{\omega}H+|h|^2_gH+Hh^2(\omega,\omega)-H^2h(\omega,\omega).$$
Since $L[f]=(\partial_t-\Delta_g)f-Hd_{\omega}f$ (for any $f$), we
see that the equation in terms of $L$ has no first-order terms:
$$L[H]=|h|^2_gH+Hh^2(\omega,\omega)-H^2h(\omega,\omega)$$

\emph{Remark.} One can also find $L[H]$ starting from the
expression:
$$L[g^{ij}h_{ij}]=L[g^{ij}]h_{ij}+g^{ij}L[h_{ij}]-2g^{kl}(\partial_kg^{ij})(\partial_lh_{ij}).$$
This may be used to check  the calculation.\newpage

\emph{Evolution of the Weingarten operator.}

The tensorial Laplacian of $S$ is the $(1,1)$ tensor $\Delta_g S$
with components $\Delta_gh^k_j$. We have:
$$\Delta_gh^k_j=g^{ik}\Delta_gh_{ij},\quad \mbox{ or }\langle
(\Delta_gS)X,Y\rangle_g=(\Delta_gh)(X,Y).$$ The evolution equation
is easily obtained:
$$(\partial_t-\Delta_g)h_j^k=(\partial_t
g^{ik})h_{ij}+g^{ik}(\partial_t-\Delta_g)h_{ij}$$
$$=H\nabla_{\omega}h_j^k+H_jh_l^k\omega^l-H_lh_j^l\omega^k+|h|^2_gh_j^k
+2HS(\omega)^kh(\omega,\partial_j)-Hh(\omega,\omega)h_j^k-Hh(S(\omega),\partial_j)\omega^k.$$
\emph{Remark:} Since the components of $\nabla S$ are given by:
$$(\nabla_{\omega}S)(\partial_j)=(\nabla_{\omega}h_j^k)\partial_k,\quad
\nabla_{\omega}h_j^k=d_{\omega}(h_j^k)+h^2(\omega,\partial_j)\omega^k-h(\omega,\partial_j)S(\omega)^k,$$
we see that upon setting $j=k$ and adding over $k$ we recover the
evolution equation for $H$.\vspace{.2cm}

The evolution equation for $h^k_j$ in terms of $L$ follows from the
calculation:
$$L[h^k_j]=L[g^{ik}]h_{ij}+g^{ik}L[h_{ij}]-2g^{mn}(\partial_mg^{ik})(\partial_nh_{ij})$$
$$=-2(\nabla_{\omega}h_m^k)h_j^m+(\partial_j|h|^2_g)\omega^k$$
$$+|h|_g^2h_j^k-Hh(\omega,\omega)h_j^k+HS(\omega)^kh(\partial_j,\omega)
+2h^3(\partial_j,\omega)\omega^k-2(h^2)^k_p\omega^ph(\partial_j,\omega).$$
Setting $j=k$ and adding over $k$, we recover the earlier expression
for $L[H]$. \vspace{.2cm}

\emph{Evolution of $|h|^2_g$.}

The fact that $g^{ij}$ is time-dependent introduces an additional
term in the usual expression:
$$(\partial_t-\Delta_g)|h|^2_g=-2|\nabla h|_g^2+2\langle
h,(\partial_t-\Delta_g)h\rangle_g+2(\partial_tg^{ij})(h^2)_{ij}.$$
Using the expressions given earlier, one easily finds:
$$(\partial_t-\Delta_g)|h|^2_g=-2|\nabla
h|_g^2+Hd_{\omega}|h|^2_g+2|h|^4_g-4Hh^3(\omega,\omega)-2H|h|^2_gh(\omega,\omega),$$
$$L[|h|^2_g]=-2|\nabla h|^2_g+2|h|^4_g-4Hh^3(\omega,\omega)-2H|h|^2_gh(\omega,\omega).$$


\begin{thebibliography}{99}

\bibitem[1]{EckerHuisken}
Ecker, K.; Huisken, G. \emph{Interior estimates for hypersurfaces
moving by mean curvature.}  Invent. Math.  \textbf{105} (1991), no.
3, 547--569.

\bibitem[2]{EckerBook}
Ecker, K. \emph{Regularity theory for mean curvature flow}. Progress
in Nonlinear Differential Equations and their Applications,
\textbf{57}. Birkhäuser Boston, Inc., Boston, MA, 2004. xiv+165 pp.
ISBN: 0-8176-3243-3

\bibitem[3]{EidelmanZhitarasu}
Eidelman, S. D.; Zhitarashu, N. V. \emph{Parabolic boundary value
problems.}  Operator Theory: Advances and Applications, 101.
Birkhäuser Verlag, Basel, 1998. xii+298 pp. ISBN: 3-7643-2972-6

\bibitem[4]{Guan} Guan, B.
\emph{Mean curvature motion of nonparametric hypersurfaces with
contact angle condition.} Elliptic and parabolic methods in geometry
(Minneapolis, MN, 1994),  47--56, A K Peters, Wellesley, MA, 1996.

\bibitem[5]{LensSeminar}
 \emph{Evolution of convex lens-shaped networks under curve shortening flow}
 O. Schn\"{u}rer, A. Azouani, M. Georgi, J. Hell, N. Jangle, A. Koeller,
 T. Marxen, S. Ritthaler, M. S\'{a}ez, F. Schulze, B. Smith,
 (Lens Seminar, FU Berlin 2007), arXiv:0711.1108

\bibitem[6]{LunardiBaconneau} Baconneau, O.; Lunardi,
A. \emph{ Smooth solutions to a class of free boundary parabolic
problems}.  Trans. Amer. Math. Soc.  \textbf{356}  (2004),  no. 3,
987--1005.

\bibitem[7]{Mantegazza et al.} Mantegazza, C.,
Novaga, M., Tortorelli, V., \emph{Motion by curvature of planar
networks.} Ann. Sc. Norm. Super. Pisa Cl. Sci. \textbf{5}(3) (2004)
no. 2, 235--324.

\bibitem[8]{MazzeoSaez} Mazzeo, R., S\'{a}ez, M. \emph{Self-similar
expanding solutions of the planar network flow}, arXiv:0704.3113

\bibitem[9]{Stahl} Stahl, A.
\emph{Convergence of solutions to the mean curvature flow with a
Neumann boundary condition.} Calc. Var. Partial Differential
Equations  \textbf{4} (1996),  no. 5, 421--441.

\bibitem[10]{Solonnikov} Solonnikov, V. A.
\emph{Lectures on evolution free boundary problems: classical
solutions}.  Mathematical aspects of evolving interfaces (Funchal,
2000),  123--175, Lecture Notes in Math., \textbf{1812}, Springer,
Berlin, 2003.

\bibitem[11]{Struwe} Struwe, M.
 \emph{The existence of surfaces of constant mean curvature with free boundaries.}
 Acta Math.  \textbf{160}  (1988),  no. 1-2, 19--64.

\end{thebibliography}
\end{document}